\documentclass[preprint,12pt]{elsarticle}

\usepackage{latexsym}
\usepackage{epsfig}
\usepackage{graphicx}
\usepackage{slashbox}
\usepackage{amsmath}
\usepackage{amssymb}
\usepackage{color}
\usepackage{amsfonts}
\usepackage{dsfont}
\usepackage{amssymb}
\newtheorem{thm}{\noindent {\bf Theorem }}
\newtheorem{lem}{\noindent {\bf Lemma }}

\newtheorem{prop}{\noindent {\bf Proposition}}
\newtheorem{example}{\noindent {\bf Example}}
\newenvironment{pro}
  {\noindent{\bf\em Proof:\ \ }}
  {\hfill$\Box$\parsep\parskip}

\journal{Elsevier}

\begin{document}


\begin{frontmatter}



\title{{\bf  Multi-symplectic Preserving Integrator for the Schr\"{o}dinger Equation with Wave Operator}}

\author[JXNU]{Linghua Kong\corref{cor1}}
\ead{konglh@mail.ustc.edu.cn}

\cortext[cor1]{Corresponding author}


\author[JXNU]{Lan Wang}

\author[LSEC]{Liying Zhang}
\ead{lyzhang@lsec.cc.ac.cn}

\address[JXNU]{School of Mathematics and Information Science, Jiangxi Normal University,\\
      Nanchang, Jiangxi, 330022, PR China}
\address[LSEC]{State Key Laboratory of Scientific and Engineering Computing, \\
    Institute of Computational Mathematics and
       Scientific/Engineering Computing, \\
    AMSS, CAS, P.O. Box 2719, Beijing, 100190, PR China}

\begin{abstract}
In the article, we discuss the conservation laws for the nonlinear
Schr\"{o}dinger equation with wave operator under multisymplectic integrator (MI).
First, the conservation laws of the continuous equation are presented and one of them is new.
The multisymplectic structure and MI are constructed for the equation. The discrete conservation laws
of the numerical method are analyzed. It is verified that the proposed MI can stably simulate the multisymplectic
Hamiltonian system excellent over long-term. It is more accurate than some energy-preserving schemes though they are of the same
accuracy. Moreover, the residual of mass is less than energy-preserving schemes under the same mesh partition over long-term.
\end{abstract}
\end{frontmatter}


\section{Introduction}
In this paper, we focus on the multisymplectic integrator (MI) for $(1+1)$ nonlinear Schr\"{o}dinger equations with wave
operator (NLSEWO) \cite{JPSJ:80}
\begin{align}\label{Schr}
  \left\{
 \begin{array}{l}
   Wu-i\alpha u_t-i\theta u_x+{\lambda}{u}+\beta{\left|u\right|^2}u=0,~(x,t)\in{[x_l, x_r]\times{(0, T]}},\\
   u(x,t)=u(x+(x_r-x_l),t),\\
   u(x,0)=f_0(x),~~ u_t(x,0)=f_1(x),~~~x\in{[x_l, x_r]},
 \end{array}\right.
\end{align}
where $Wu=u_{tt}-u_{xx}+\gamma{u_{tx}}$, $i^2=-1$, $\alpha, \gamma, \theta, \lambda$ and $\beta$ are real constants,
at the same time, $f_0(x)$ and $f_1(x)$ are given functions. The equation describes the nonlinear interaction between
two quasi-monochromatic waves. It is one of the non-resonant interaction.
\begin{prop}
The determined problem (\ref{Schr})
satisfies the following conservation laws
\begin{itemize}
 \item Energy invariant
  \begin{align}\label{energy}
    \newcommand\diff{\,{\mathrm d}}
      \mathcal{E}(t)=\int_{x_l}^{x_r}( \left|u_{t}\right|^2+\left|u_{x}\right|^2+i\theta{u}\bar{u}_{x}+\lambda\left|u\right|^2
        +\frac{\beta}{2}\left|u\right|^4)\diff x=\mathcal{E}(0).
  \end{align}

\item Mass invariant
\begin{align}\label{mass}
  \mathcal{Q}(t)=\int_{x_l}^{x_r}\left[(u_t\overline{u}-\overline{u_t}u)-\gamma u\overline{u_x}-i\alpha|u|^2\right]dx=\mathcal{Q}(0).
\end{align}
\end{itemize}

\end{prop}
\begin{pro}
We multiply Eq. (\ref{Schr}) with $\overline{u_t}$ and integrate
over spatial domain
\begin{align}\label{Eq2}
 \int_{x_l}^{x_r}(u_{tt}\overline{u_t}-u_{xx}\overline{u_t}+\gamma{u_{tx}}\overline{u_t}-
 i\alpha{u_t}\overline{u_t}-i\theta{u_x}\overline{u_t}+{\lambda}{u}\overline{u_t}+\beta{\left|u\right|^2}u\overline{u_t})dx=0.
\end{align}
Next, we multiply the conjugation of Eq. (\ref{Schr}) with ${u_t}$
\begin{align}\label{Eq3}
 \int_{x_l}^{x_r}(\overline{u_{tt}}{u_t}-\overline{u_{xx}}{u_t}+\gamma\overline{u_{tx}}{u_t}+
 i\alpha\overline{u_t}{u_t}+i\theta\overline{u_x}{u_t}+{\lambda}\overline{u}{u_t}+\beta{\left|u\right|^2}\overline{u}{u_t})dx=0.
\end{align}
 Adding Eq. (\ref{Eq2}) to Eq. (\ref{Eq3}) with integration by part under the boundary conditions, we have
\begin{align*}
  \frac{d}{dt}\int_{x_l}^{x_r}\left[\left|u_{t}\right|^2+\left|u_{x}\right|^2
   +i\theta({u}\overline{u_x})+{\lambda}\left|u\right|^2+\frac{\beta}{2}{\left|u\right|^4}\right]dx=0.
\end{align*}
This is just what we desire.
Consequently, we prove another invariant:
\begin{align*}
    & \int_{x_l}^{x_r}\left[(u_{tt}\overline{u}-u_{xx}\overline{u}+\gamma u_{xt}\overline{u}-i\alpha u_t\overline{u}
        -i\theta u_x\overline{u}+\lambda |u|^2+\beta |u|^4)\right.\\
    &    \qquad -\left.(\overline{u_{tt}}u-\overline{u_{xx}}u+\gamma\overline{u_{xt}}u+i\alpha \overline{u_t}u+
       i\theta \overline{u_x}u+\lambda |u|^2+\beta |u|^4)\right]dx\\
  = & \frac{d}{dt}\int_{x_l}^{x_r}\left[(u_t\overline{u}-\overline{u_t}u)-\gamma (u\overline{u_x})-i\alpha |u|^2\right]\\
  = & 0.
\end{align*}
The proof is completed.
\end{pro}

It is noted that the mass invariant (\ref{mass}) is {\bf new} which did not appear in existing literatures to our knowledge.

In \cite{GBL}, Guo proposed an implicit nonconservative difference scheme for NLSEWO.
Based on the first conservative quantity, Zhang et al developed some energy-preserving schemes for it
\cite{ZLM:AMC:03, WTC:AMC:06, WSS:JCAM:11} in case of $\gamma=\theta=\lambda=0$. Wang considered its Fourier
pseudo-spectral method under multisymplectic context \cite{WJ:JCM:07}.

At the end of last century, MIs have been put forward and application to large numbers of
partial differential equations, such as wave equation \cite{ Bridge:PLA:01, Reich:JCP:00,WYS:CS:02}, nonlinear
Schr\"{o}dinger-type equations \cite{HJL:CICP:10, HJL:ANM:06},Dirac equation \cite{HJL:JCP:06},
Maxwell's equations \cite{KLH:JCP:10}, RLW equation \cite{CJX:CPC:10}. The most important character of MIs is its multisymplecticity, and
other conservative properties are preserved excellently despite of not exactly \cite{H:JCP:07, HJL:JCP:06, Schober:06, JPA:06, ZHJ:ANM:11}.
In the article, we investigate the MIs and its global conservative properties.

The rest of the paper is organized as follows: In Section \ref{Preliminary}, some preliminary knowledge is prepared for which will often be
used later.  In Section \ref{method}, it presents the multisymplectic structure and an MI for
NLSEWO. The conservation properties of the proposed numerical schemes are investigated in
Section \ref{analysis}. In Section \ref{example}, we present some
numerical examples and detailed numerical results. Some conclusions
are given to end this paper.

\section{Preliminary Knowledge}\label{Preliminary}
In this section, we give some notations and knowledge we will frequently be used.
  A uniform partition of the domain under consideration is
$$x_k=x_l+kh, t_j=j\tau, k=0, 1, 2, \cdots, K; j=0, 1, \cdots, J.$$
Where $h=\frac{x_r-x_l}{K}$ and $\tau=\frac{T}{J}$ denote the spatial mesh size and temporal step length.
$u_k^j$ is an approximation of $u(x,t)$ at the node $(x_k, t_j)$, and $U_k^j=u(x_k, t_j)$. Some notations about
difference quotient:
\begin{align*}
    & {u_k^j}_x=\frac{u_{k+1}^j-u_k^j}{h},  {u_k^j}_{\bar{x}}=\frac{u_{k}^j-u_{k-1}^j}{h},
      {u_{k+\frac{1}{2}}^{j}}_{\hat{x}}=\frac{u_{k+1}^j-u_{k}^j}{h}, {u_k^j}_{2x}=\frac{u_{k+1}^j-u_{k-1}^j}{2h};\\
    & {u_k^j}_t=\frac{u_{k}^{j+1}-u_k^j}{\tau},  {u_k^j}_{\bar{t}}=\frac{u_{k}^j-u_{k}^{j-1}}{\tau},
      {u_k^{j+\frac{1}{2}}}_{\hat{t}}=\frac{u_{k}^{j+1}-u_{k}^{j}}{\tau},  {u_k^j}_{2t}=\frac{u_{k}^{j+1}-u_{k}^{j-1}}{2\tau};\\
    & \delta_t^2u_k^j=\frac{u_{k}^{j+1}-2u_{k}^{j}+u_{k}^{j-1}}{\tau^2}, \delta_x^2u_k^j=\frac{u_{k+1}^{j}-2u_{k}^{j}+u_{k-1}^{j}}{h^2}.
\end{align*}

For any vectors $u^j, v^j\in \mathbb{C}^K$, the inner product and norms are defined as:
\begin{align*}
  & \langle u^j, v^j\rangle=h\sum\limits_ku_k^j\overline{v_k^j}, ~~~\left\|u^j\right\|^2=\langle u^j, u^j\rangle,\\
  & \left\|u^j\right\|_{\frac{1}{2}}^2=h\sum\limits_k\left|u_{k+\frac{1}{2}}^j\right|^2, ~~~
    \left\|u^j\right\|_\infty=\max\limits_k\left|u_k^j\right|.
\end{align*}

\begin{lem}\label{Green}
  For all complex mesh functions $\{u_{k}^{j}\}$ and $\{v_{k}^{j}\}$ with periodic or homogeneous boundary condition, we have the following conclusions:
  \begin{itemize}
   \item Discrete Green formula:
    \begin{align}
      \langle \delta^2_{x}u^j, v^j\rangle=-\langle \delta_{\bar{x}}u^j, \delta_{\bar{x}}v^j\rangle;
    \end{align}
   \item $\langle \delta_{\bar{x}}u^j, u^j\rangle=-\langle u^j, \delta_{\bar{x}}u^j\rangle$ is purely imaginary.
   \item $Re\langle \delta_t^2u^j, u_{2t}^j\rangle=\frac{1}{2}\left\|u_t^j\right\|^2_{\bar{t}}$, where `$Re$' denotes taking real part.
  \end{itemize}
\end{lem}
These conclusion can be verified easily.

\section{Construction of Multi-symplectic scheme}\label{method}
In this section, we firstly describe the multi-symplectic structure
and local conservation laws for the Eq.(\ref{Schr}). In order to
rewrite the complex equations as a real one, we suppose that
\begin{align*}
  \left\{
    \begin{array}{l}
      u(x,t)=\varphi(x,t)+i\psi(x,t),\\
      u_t(x,t)=v(x,t)+iw(x,t),\\
      u_x(x,t)=f(x,t)+ig(x,t),
   \end{array}
  \right.
\end{align*}
 where $\varphi(x,t), \psi(x,t), v(x,t), w(x,t),
f(x,t), g(x,t)$ are all real-valued functions.
 Let $z=(\varphi,\psi,v,w,f,g)^T$, we have multi-symplectic equations
 \begin{align}\left\{\label{sche}
 \begin{array}{ll}
    \alpha{\psi_t}+v_t+\frac{\gamma}{2}f_t+\theta{\psi_x}+\frac{\gamma}{2}v_x-f_x=-{\lambda}\varphi-\beta(\varphi^2+\psi^2)\varphi,\\
    -\alpha{\varphi_t}+w_t+\frac{\gamma}{2}g_t-\theta{\varphi_x}+\frac{\gamma}{2}w_x-g_x=-{\lambda}\psi-\beta(\varphi^2+\psi^2)\psi,\\
    -\varphi_t-\frac{\gamma}{2}\varphi_x=-v-\frac{\gamma}{2}f,\\
    -\psi_t-\frac{\gamma}{2}\psi_x=-w-\frac{\gamma}{2}g,\\
    -\frac{\gamma}{2}\varphi_t+\varphi_x=-\frac{\gamma}{2}v+f,\\
    -\frac{\gamma}{2}\psi_t+\psi_x=-\frac{\gamma}{2}w+g.
 \end{array}\right.
\end{align}
Then, we can cast (\ref{sche}) into the multi-symplectic framework
\begin{align}\label{mtlaw}
Mz_t+Kz_x=\bigtriangledown _zS(z),
\end{align}
where $\bigtriangledown$ is the gradient operator.
The Hamiltonian function is
\begin{align*}
 S(z)=-\frac{1}{2}[{{\lambda}}(\varphi^2+\psi^2)+\frac{\beta}{2}(\varphi^2+\psi^2)^2+v^2+w^2-(f^2+g^2)+\gamma(vf+wg)].
\end{align*}
\begin{align*}
M=\left[
    \begin{array}{cccccc}
      0          &  \alpha     &¡¡1     &  0  &\frac{1}{2}\gamma &   0   \\
    -\alpha      &  0          &  0     &  1    & 0              & \frac{1}{2}\gamma \\
     -1          &  0          &  0     &  0    & 0              &   0  \\
      0          &  -1         &  0     &  0    & 0              &   0\\
-\frac{\gamma}{2}&   0         &  0     &  0    & 0              &   0\\
      0      &-\frac{\gamma}{2}&  0     &  0    & 0              &   0\\
  \end{array}
  \right], ~
K=\left[
    \begin{array}{cccccc}
     0            &  \theta           &\frac{\gamma}{2} &  0                    &-1    & 0   \\
-\theta           &  0                & 0               & \frac{\gamma }{2}     & 0    & -1  \\
-\frac{\gamma }{2}&  0                & 0               &  0                    & 0    & 0   \\
     0            &-\frac{\gamma }{2}  &0                &  0                    & 0    & 0   \\
     1            &  0                & 0               &  0                    & 0    & 0   \\
     0            &  1               &0                &  0                    & 0    & 0   \\
  \end{array}
  \right].
\end{align*}

According to the multi-symplectic theoretical background, the
multi-symplectic system (\ref{mtlaw}) satisfies local conservation
laws as follows:
\begin{itemize}
 \item Multi-symplectic conservation law
 \begin{align}
   \frac{\partial}{\partial{t}}\omega+\frac{\partial}{\partial{x}}\kappa=0, ~~ \forall (x, t),
\end{align}
where $\omega $ and $\kappa$ are pre-symplectic 2-forms
\begin{align*}
  &¡¡\omega={\alpha}d\varphi\wedge{d\psi}+d\varphi\wedge{dv}+\frac{\gamma}{2}d\varphi\wedge{df}+d\psi\wedge{dw}+\frac{\gamma}{2}d\psi\wedge{dg},\\
¡¡&  \kappa={\theta}d\varphi\wedge{d\psi}-d\varphi\wedge{df}+\frac{\gamma}{2}d\varphi\wedge{dv}-d\psi\wedge{dg}+\frac{\gamma}{2}d\psi\wedge{dw}.
\end{align*}

\item Local energy conservation law:
\begin{align}
\frac{\partial}{\partial{t}}E(z)+\frac{\partial}{\partial{x}}F(z)=0, ~~  \forall (x, t),
\end{align}
where the energy density $E(z)$ and the energy flux $F(z)$ are
\begin{align*}
  & E(z)=\frac{1}{2}[{\lambda}(\varphi^2+\psi^2)+\frac{\beta}{2}(\varphi^2+\psi^2)^2+v^2+w^2+f^2+g^2+\theta(\varphi g-\psi f)],\\
  & F(z)=\frac{\theta}{2}(\varphi w-\psi v)-\frac{\gamma}{2}(v^2+w^2).
\end{align*}

\item Local momentum conservation law:
\begin{align}
\frac{\partial}{\partial{t}}I(z)+\frac{\partial}{\partial{x}}G(z)=0,
\end{align}
where the momentum density $I(z)$ and the momentum flux $G(z)$ are
\begin{align*}
  & I(z)=\frac{\alpha}{2}(\varphi g-\psi f)-(fv+gw)-\frac{\gamma}{2}(f^2+g^2),\\
  & G(z)=-\frac{\lambda}{2}(\varphi^2+\psi^2)-\frac{\beta}{4}(\varphi^2+\psi^2)^2+\frac{1}{2}(v^2+w^2+f^2+g^2)-\frac{\alpha}{2}(\varphi w-\psi v).~~~~~~~~~~~~~~~~~~~~~~~~~~
\end{align*}
\end{itemize}

It is well known that the local conservation laws imply that the
density can be various, but the increment of the density in time
puts up with the flux in space.

Applying the multisymplectic midpoint integrator
\begin{equation}\label{box}
     M\delta_{\hat{t}z_{k+\frac{1}{2}}}^{j+\frac{1}{2}}+K\delta_{\hat{x}}z_{k+\frac{1}{2}}^{j+\frac{1}{2}}
    =\nabla_zS(z_{k+\frac{1}{2}}^{j+\frac{1}{2}}),
\end{equation}
where $z_{k+\frac{1}{2}}^{j+\frac{1}{2}}=\frac{1}{2}(z_{k+1}^{j+\frac{1}{2}}+z_{k}^{j+\frac{1}{2}})=
\frac{1}{2}(z_{k+\frac{1}{2}}^{j+1}+z_{k+\frac{1}{2}}^{j})
=\frac{1}{4}(z_{k+1}^{j+1}+z_{k+1}^{j}+z_{k}^{j+1}+z_{k}^{j})$, to the multisymplectic Hamiltonian system (\ref{sche}), one has
\begin{align}\label{midpoint}
 \left\{
 \begin{array}{rl}
    & \alpha\delta_{\hat{t}}{\psi_{k+\frac{1}{2}}^{j+\frac{1}{2}}}+\delta_{\hat{t}}{v_{k+\frac{1}{2}}^{j+\frac{1}{2}}}
      +\frac{\gamma}{2}\delta_{\hat{t}}{f_{k+\frac{1}{2}}^{j+\frac{1}{2}}}
      +\theta{\delta_{\hat{x}}{\psi_{k+\frac{1}{2}}^{j+\frac{1}{2}}}}
      +\frac{\gamma}{2}\delta_{\hat{x}}{v_{k+\frac{1}{2}}^{j+\frac{1}{2}}}
      -\delta_{\hat{x}}{f_{k+\frac{1}{2}}^{j+\frac{1}{2}}}\\
  = & -{\lambda}{\varphi_{k+\frac{1}{2}}^{j+\frac{1}{2}}}
      -\beta\left[\left(\varphi_{k+\frac{1}{2}}^{j+\frac{1}{2}}\right)^2+\left(\psi_{k+\frac{1}{2}}^{j+\frac{1}{2}}\right)^2\right]{\varphi_{k+\frac{1}{2}}^{j+\frac{1}{2}}},\\
    & -\alpha\delta_{\hat{t}}{\varphi_{k+\frac{1}{2}}^{j+\frac{1}{2}}}
      +\delta_{\hat{t}}{w_{k+\frac{1}{2}}^{j+\frac{1}{2}}}
      +\frac{\gamma}{2}\delta_{\hat{t}}{g_{k+\frac{1}{2}}^{j+\frac{1}{2}}}
      -\theta\delta_{\hat{x}}{\varphi_{k+\frac{1}{2}}^{j+\frac{1}{2}}}
      +\frac{\gamma}{2}\delta_{\hat{x}}{w_{k+\frac{1}{2}}^{j+\frac{1}{2}}}
      -\delta_{\hat{x}}{g_{k+\frac{1}{2}}^{j+\frac{1}{2}}}\\
  = & -{\lambda}{\psi_{k+\frac{1}{2}}^{j+\frac{1}{2}}}
      -\beta\left[\left({\varphi_{k+\frac{1}{2}}^{j+\frac{1}{2}}}\right)^2+\left({\psi_{k+\frac{1}{2}}^{j+\frac{1}{2}}}\right)^2\right]{\psi_{k+\frac{1}{2}}^{j+\frac{1}{2}}},\\
    & -\delta_{\hat{t}}{\varphi_{k+\frac{1}{2}}^{j+\frac{1}{2}}}=-{v_{k+\frac{1}{2}}^{j+\frac{1}{2}}}, ~~
      -\delta_{\hat{t}}{\psi_{k+\frac{1}{2}}^{j+\frac{1}{2}}}=-{w_{k+\frac{1}{2}}^{j+\frac{1}{2}}},\\
    & \delta_{\hat{x}}{\varphi_{k+\frac{1}{2}}^{j+\frac{1}{2}}}={f_{k+\frac{1}{2}}^{j+\frac{1}{2}}}, ~~
      \delta_{\hat{x}}{\psi_{k+\frac{1}{2}}^{j+\frac{1}{2}}}={g_{k+\frac{1}{2}}^{j+\frac{1}{2}}}.
 \end{array}\right.
\end{align}
For the details of the method and the theoretical results on local conservation laws of the numerical method, we refer to \cite{HJL:JCP:06} and
references therein.

The MI (\ref{midpoint}) is of second order both in time and space. After tedious calculation,
by eliminating the introduced variables, it can be reformulated into
\begin{align}\label{Preissman}
 \begin{array}{rl}
     & \frac{1}{2}(\delta_t^2{u_{k+\frac{1}{2}}^{j}}+\delta_t^2u_{k-\frac{1}{2}}^{j})-
       \frac{1}{2}(\delta_x^2{u_{k}^{j+\frac{1}{2}}}+\delta_x^2u_{k}^{j-\frac{1}{2}})
       -\frac{i\alpha}{2}(\delta_{2t}{u_{k+\frac{1}{2}}^{j}}+\delta_{2t}u_{k-\frac{1}{2}}^{j})\\
     & -\frac{i\theta}{2}(\delta_{2x}{u_{k}^{j+\frac{1}{2}}}+\delta_{2x}u_{k}^{j-\frac{1}{2}})
       +{\gamma}{\delta_{2t}\delta_{2x}u_k^j}
       +\frac{{\lambda}}{4}\left[u_{k+\frac{1}{2}}^{j+\frac{1}{2}}+u_{k-\frac{1}{2}}^{j+\frac{1}{2}}
       +u_{k+\frac{1}{2}}^{j-\frac{1}{2}}+u_{k-\frac{1}{2}}^{j-\frac{1}{2}}\right]\\
    & +\frac{\beta}{4}\left[\left|u_{k+\frac{1}{2}}^{j+\frac{1}{2}}\right|^2u_{k+\frac{1}{2}}^{j+\frac{1}{2}}
       +\left|u_{k+\frac{1}{2}}^{j-\frac{1}{2}}\right|^2u_{k+\frac{1}{2}}^{j-\frac{1}{2}}+\left|u_{k-\frac{1}{2}}^{j+\frac{1}{2}}\right|^2u_{k-\frac{1}{2}}^{j+\frac{1}{2}}
       +\left|u_{k-\frac{1}{2}}^{j-\frac{1}{2}}\right|^2u_{k-\frac{1}{2}}^{j-\frac{1}{2}}\right]=0.
 \end{array}
\end{align}
By Taylor expansion, the MI (\ref{Preissman}) is of second order both in space and time, that is, the truncation error is
$\mathcal{T}_k^{j+\frac{1}{2}}=\mathcal{O}(\tau^2+h^2)$.

\section{Conservation Laws Analysis}\label{analysis}
In this section, the theoretical analysis about the MI (\ref{Preissman}) is derived. It is suggested that the MI can preserve the energy and mass
very well though they are not exactly.

 Firstly, we investigate the discrete energy
conservation law. To the purpose, we multiply the Eq. (\ref{Preissman}) with
$$2\left(\overline{u_{k}^{j+\frac{1}{2}}}-\overline{u_{k}^{j-\frac{1}{2}}}\right)
=\overline{u_k^{j+1}}-\overline{u_k^{j-1}}=2\tau\delta_{2t}{\overline{u_k^j}}
=\tau(\delta_{\hat{t}}\overline{u_{k}^{j+\frac{1}{2}}}+\delta_{\hat{t}}\overline{u_{k}^{j-\frac{1}{2}}}) ,$$ and
sum over index $k$. Then, from the third formulate in Lemma \ref{Green}, one can obtain the real part of the first term is
\begin{align}\label{energy_first}
  \left\|\delta_{\hat{t}}u^{j+\frac{1}{2}}\right\|_\frac{1}{2}^2-\left\|\delta_{\hat{t}}u^{j-\frac{1}{2}}\right\|_\frac{1}{2}^2.
\end{align}
By Green formula, the real part of the second term is
\begin{align}\label{energy_second}
 \begin{array}{rl}
     & -\frac{1}{2}Re\left\langle\delta_x^2{u^{j+\frac{1}{2}}}+\delta_x^2u^{j-\frac{1}{2}}, u^{j+\frac{1}{2}}-u_k^{j-\frac{1}{2}}\right\rangle\\
  =  & Re\left\langle\delta_{\bar{x}}u^{j+\frac{1}{2}}+\delta_{\bar{x}}u^{j-\frac{1}{2}},
       \delta_{\bar{x}}u^{j+\frac{1}{2}}-\delta_{\bar{x}}u^{j-\frac{1}{2}}\right\rangle\\
  =  & \left\|\delta_{\bar{x}}u^{j+\frac{1}{2}}\right\|^2-\left\|\delta_{\bar{x}}u^{j-\frac{1}{2}}\right\|^2.
 \end{array}
\end{align}

The third term
\begin{align*}
 \begin{array}{rl}
     -\frac{i\alpha}{2}h\sum\limits_k\left[\delta_{2t}{u_{k+\frac{1}{2}}^{j}}+\delta_{2t}u_{k-\frac{1}{2}}^{j}\right]2\tau\delta_{2t}{\overline{u_k^j}}
   = & -2\tau{i\alpha}h\sum\limits_k\left|\delta_{2t}{u_{k+\frac{1}{2}}^{j}}\right|^2\\
   = & -2\tau{i\alpha}\left\|\delta_{2t}{u^{j}}\right\|_\frac{1}{2}^2,
 \end{array}
\end{align*}
is purely imaginary. Judged from the second conclusion of Lemma \ref{Green}, the  fifth term
${\gamma}h\sum\limits_k(\delta_{2x}\delta_{2t}u_k^j)2\tau\delta_{2t}\overline{u_k^j}$
is purely imaginary, too. Now, we analyze the fourth term
\begin{align*}
   &  -\frac{i\theta}{2}h\sum\limits_k(\delta_{2x}{u_{k}^{j+\frac{1}{2}}}+\delta_{2x}u_{k}^{j-\frac{1}{2}})(\overline{u_k^{j+1}}-\overline{u_k^{j-1}})\\
 = & -{i\theta}h\sum\limits_k(\delta_{2x}{u_{k}^{j+\frac{1}{2}}}+\delta_{2x}u_{k}^{j-\frac{1}{2}})
     (\overline{u_k^{j+\frac{1}{2}}}-\overline{u_k^{j-\frac{1}{2}}})\\
 = & -\frac{i\theta}{2{h}}h\sum\limits_k\left[(u_{k+1}^{j+\frac{1}{2}}+u_{k+1}^{j-\frac{1}{2}})
     (\overline{u_k^{j+\frac{1}{2}}}-\overline{u_k^{j-\frac{1}{2}}})
     -(u_{k}^{j+\frac{1}{2}}+u_{k}^{j-\frac{1}{2}})(\overline{u_{k+1}^{j+\frac{1}{2}}}-\overline{u_{k+1}^{j-\frac{1}{2}}})\right].
\end{align*}
 Its real part is
 \begin{align}\label{energy_fourth}
  \begin{array}{rl}
    &  -\frac{i\theta}{2h}h\sum\limits_k\left[(u_{k+1}^{j+\frac{1}{2}}\overline{u_k^{j+\frac{1}{2}}}-u_{k}^{j+\frac{1}{2}}\overline{u_{k+1}^{j+\frac{1}{2}}})
       -(u_{k+1}^{j-\frac{1}{2}}\overline{u_k^{j-\frac{1}{2}}}-u_{k}^{j-\frac{1}{2}}\overline{u_{k+1}^{j-\frac{1}{2}}})\right]\\
  = &  -\frac{i\theta}{2h}h\sum\limits_k\left[(u_{k+1}^{j+\frac{1}{2}}\overline{u_k^{j+\frac{1}{2}}}
       -u_{k+1}^{j+\frac{1}{2}}\overline{u_{k+1}^{j+\frac{1}{2}}})
       -(u_{k}^{j+\frac{1}{2}}\overline{u_{k+1}^{j+\frac{1}{2}}}-u_{k}^{j+\frac{1}{2}}\overline{u_{k}^{j+\frac{1}{2}}})\right.\\
    &  \qquad-(u_{k+1}^{j-\frac{1}{2}}\overline{u_k^{j-\frac{1}{2}}}-u_{k+1}^{j-\frac{1}{2}}\overline{u_{k+1}^{j-\frac{1}{2}}})
       \left.+(u_{k}^{j-\frac{1}{2}}\overline{u_{k+1}^{j-\frac{1}{2}}}-u_{k}^{j-\frac{1}{2}}\overline{u_{k}^{j-\frac{1}{2}}})\right]\\
  = & i\theta h\sum\limits_k\left[u_{k+\frac{1}{2}}^{j+\frac{1}{2}}\delta_{\hat{x}}\overline{u_{k+\frac{1}{2}}^{j+\frac{1}{2}}}
      -u_{k+\frac{1}{2}}^{j-\frac{1}{2}}\delta_{\hat{x}}\overline{u_{k+\frac{1}{2}}^{j-\frac{1}{2}}}\right].
 \end{array}
\end{align}
The real part of the sixth term reads
\begin{align}\label{energy_sixth}
 \begin{array}{rl}
   & \frac{{2\lambda}}{4}h Re\sum\limits_k\left[\left(u_{k+\frac{1}{2}}^{j+\frac{1}{2}}+u_{k+\frac{1}{2}}^{j-\frac{1}{2}}\right)
     +\left(u_{k-\frac{1}{2}}^{j+\frac{1}{2}}+u_{k-\frac{1}{2}}^{j-\frac{1}{2}}\right)\right]
       \left(\overline{u_k^{j+\frac{1}{2}}}-\overline{u_k^{j-\frac{1}{2}}}\right)\\
 = & {\lambda}h Re\sum\limits_k\left(u_{k+\frac{1}{2}}^{j+\frac{1}{2}}+u_{k+\frac{1}{2}}^{j-\frac{1}{2}}\right)
     \left(\overline{u_{k+\frac{1}{2}}^{j+\frac{1}{2}}}-\overline{u_{k+\frac{1}{2}}^{j-\frac{1}{2}}}\right)\\
 = & \lambda\left(\left\|u^{j+\frac{1}{2}}\right\|_\frac{1}{2}^2-\left\|u^{j-\frac{1}{2}}\right\|_\frac{1}{2}^2\right).
 \end{array}
\end{align}

The last term is
\begin{align*}
   & \frac{2\beta}{{4}}h\sum\limits_k\left[\left|u_{k+\frac{1}{2}}^{j+\frac{1}{2}}\right|^2
      u_{k+\frac{1}{2}}^{j+\frac{1}{2}}+\left|u_{k+\frac{1}{2}}^{j-\frac{1}{2}}\right|^2u_{k+\frac{1}{2}}^{j-\frac{1}{2}}
      +\left|u_{k-\frac{1}{2}}^{j+\frac{1}{2}}\right|^2u_{k-\frac{1}{2}}^{j+\frac{1}{2}}
      +\left|u_{k-\frac{1}{2}}^{j-\frac{1}{2}}\right|^2u_{k-\frac{1}{2}}^{j-\frac{1}{2}}\right]\left[\overline{u_k^{j+\frac{1}{2}}}-\overline{u_k^{j-\frac{1}{2}}}\right]\\
 = &  \beta h\sum\limits_k\left[\left|u_{k+\frac{1}{2}}^{j+\frac{1}{2}}\right|^2
      u_{k+\frac{1}{2}}^{j+\frac{1}{2}}+\left|u_{k+\frac{1}{2}}^{j-\frac{1}{2}}\right|^2u_{k+\frac{1}{2}}^{j-\frac{1}{2}}\right]
      \left[\overline{u_{k+\frac{1}{2}}^{j+\frac{1}{2}}}-\overline{u_{k+\frac{1}{2}}^{j-\frac{1}{2}}}\right]\\
 = &  {\beta}h\sum\limits_k\left[\left|u_{k+\frac{1}{2}}^{j+\frac{1}{2}}\right|^4-\left|u_{k+\frac{1}{2}}^{j-\frac{1}{2}}\right|^4
      +\left|u_{k+\frac{1}{2}}^{j-\frac{1}{2}}\right|^2u_{k+\frac{1}{2}}^{j-\frac{1}{2}}\overline{u_{k+\frac{1}{2}}^{j+\frac{1}{2}}}
      -\left|u_{k+\frac{1}{2}}^{j+\frac{1}{2}}\right|^2u_{k+\frac{1}{2}}^{j+\frac{1}{2}}\overline{u_{k+\frac{1}{2}}^{j-\frac{1}{2}}}\right],
\end{align*}
with the real part
\begin{align*}
   & 2\frac{\beta}{2}\left[\left\|u^{j+\frac{1}{2}}\right\|_\frac{1}{2}^4-\left\|u^{j-\frac{1}{2}}\right\|_\frac{1}{2}^4\right]
     +\beta h Re \sum\limits_k\left(\left|u_{k+\frac{1}{2}}^{j-\frac{1}{2}}\right|^2
      u_{k+\frac{1}{2}}^{j-\frac{1}{2}}\overline{u_{k+\frac{1}{2}}^{j+\frac{1}{2}}}
      -\left|u_{k+\frac{1}{2}}^{j+\frac{1}{2}}\right|^2u_{k+\frac{1}{2}}^{j+\frac{1}{2}}\overline{u_{k+\frac{1}{2}}^{j-\frac{1}{2}}}\right)\\
 = & \frac{\beta}{2}\left[\left\|u^{j+\frac{1}{2}}\right\|_\frac{1}{2}^4-\left\|u^{j-\frac{1}{2}}\right\|_\frac{1}{2}^4\right]+
     \frac{\beta}{2}h\sum\limits_k\left[\left|u_{k+\frac{1}{2}}^{j+\frac{1}{2}}\right|^2-\left|u_{k+\frac{1}{2}}^{j-\frac{1}{2}}\right|^2\right]
     \left|u_{k+\frac{1}{2}}^{j+\frac{1}{2}}-u_{k+\frac{1}{2}}^{j-\frac{1}{2}}\right|^2.
\end{align*}

In summary, one has the following theorem:
\begin{thm}\label{discrete_energy}
  The Preissman MI (\ref{Preissman}) possesses the discrete implicit
  energy conservation law, i.e.
\begin{align}
      \mathcal{E}^{j+\frac{1}{2}}-\mathcal{E}^{j-\frac{1}{2}}
   =-\frac{\beta}{2}h\sum\limits_k\left[\left|u_{k+\frac{1}{2}}^{j+\frac{1}{2}}\right|^2-\left|u_{k+\frac{1}{2}}^{j-\frac{1}{2}}\right|^2\right]
     \left|\delta_{2t}u_{k+\frac{1}{2}}^j\right|^2\tau^2,
\end{align}
where
$
   \mathcal{E}^{j+\frac{1}{2}}=\|\delta_{\hat{t}}u^{j+\frac{1}{2}}\|_\frac{1}{2}^2
    +i{\theta}h\sum\limits_ku_{k+\frac{1}{2}}^{j+\frac{1}{2}}\delta_{\hat{x}}\overline{u_{k+\frac{1}{2}}^{j+\frac{1}{2}}}
    +\|\delta_{\hat{x}}{u^{j+\frac{1}{2}}}\|_\frac{1}{2}^2+{\lambda}\left\|u^{j+\frac{1}{2}}\right\|_\frac{1}{2}^2
    +\frac{\beta}{2}\left\|u^{j+\frac{1}{2}}\right\|_\frac{1}{2}^4$.
In particular, in case of $\beta=0$, the conservation law is explicit, that is,
\begin{align}
  \mathcal{E}^{j+\frac{1}{2}}=\mathcal{E}^{j-\frac{1}{2}}=\cdots=\mathcal{E}^{\frac{1}{2}}.
\end{align}
\end{thm}
Furthermore, we have the implicit mass conservation law:
\begin{thm}\label{discrete_mass}
  The Preissman MI (\ref{Preissman}) admits the  implicit mass conservation law
\begin{align}\label{implicit_mass}
 \begin{array}{rl}
     & \frac{\mathcal{Q}^{j+\frac{1}{2}}-\mathcal{Q}^{j-\frac{1}{2}}}{\tau} \\
  =  & -\frac{\beta}{2}h\sum\limits_k\left[\left|u_{k+\frac{1}{2}}^{j+\frac{1}{2}}\right|^2-\left|u_{k+\frac{1}{2}}^{j-\frac{1}{2}}\right|^2\right]
       \left(u_{k+\frac{1}{2}}^{j+\frac{1}{2}}-u_{k+\frac{1}{2}}^{j-\frac{1}{2}}\right)
       \left(\overline{u_{k+\frac{1}{2}}^{j+\frac{1}{2}}}+\overline{u_{k+\frac{1}{2}}^{j-\frac{1}{2}}}\right)\\
     & +\frac{\beta}{2}h\sum\limits_k\left[\left|u_{k+\frac{1}{2}}^{j+\frac{1}{2}}\right|^2-\left|u_{k+\frac{1}{2}}^{j-\frac{1}{2}}\right|^2\right]^2,
 \end{array}
\end{align}
  where
$\mathcal{Q}^{j+\frac{1}{2}}=h\sum\limits_k\left[\delta_{\hat{t}}u_{k+\frac{1}{2}}^{j+\frac{1}{2}}\overline{u_{k+\frac{1}{2}}^{j+\frac{1}{2}}}-
  u_{k+\frac{1}{2}}^{j+\frac{1}{2}}\delta_{\hat{t}}\overline{u_{k+\frac{1}{2}}^{j+\frac{1}{2}}}\right]-
  \gamma h\sum\limits_ku_{k+\frac{1}{2}}^{j+\frac{1}{2}}\delta_{\hat{x}}\overline{u_{k+\frac{1}{2}}^{j+\frac{1}{2}}}
  -i\alpha\left\|u^{j+\frac{1}{2}}\right\|^2$.
In particular, if $\beta=0$, the conservation law is explicit, that is,
\begin{align}
  \mathcal{Q}^{j+\frac{1}{2}}=\mathcal{Q}^{j-\frac{1}{2}}=\cdots=\mathcal{Q}^{\frac{1}{2}}.
\end{align}
\end{thm}
\begin{pro}
 Taking inner product of (\ref{Preissman}) with
  $$B{u_k^j}={u_k^{j+1}}+2{u_k^j}+{u_k^{j-1}}
   =2({u_{k}^{j+\frac{1}{2}}}+{u_{k}^{j-\frac{1}{2}}}),$$
the first term is showed as
 \begin{align*}
    &  \frac{1}{2}h\sum\limits_k(\delta_t^2{u_{k+\frac{1}{2}}^{j}}+\delta_t^2u_{k-\frac{1}{2}}^{j})
      (\overline{u_k^{j+1}}+2\overline{u}_k^j+\overline{u_k^{j-1}})\\
  = & \frac{1}{\tau^2}h\sum\limits_k\left[(u_{k+\frac{1}{2}}^{j+1}-u_{k+\frac{1}{2}}^{j})-(u_{k+\frac{1}{2}}^{j}-u_{k+\frac{1}{2}}^{j-1})\right]
       \left[(\overline{u_{k+\frac{1}{2}}^{j+1}}+\overline{u_{k+\frac{1}{2}}^{j}})+
       (\overline{u_{k+\frac{1}{2}}^{j}}+\overline{u_{k+\frac{1}{2}}^{j-1}})\right],
 \end{align*}
 whose imaginary part is as follows:
\begin{align}\label{mass_first}
 \begin{array}{rl}
    &  \frac{2}{\tau^2}h\sum\limits_k\left[u_{k+\frac{1}{2}}^{j+1}\overline{u_{k+\frac{1}{2}}^{j}}-\overline{u_{k+\frac{1}{2}}^{j+1}}u_{k+\frac{1}{2}}^{j}
       -u_{k+\frac{1}{2}}^{j}\overline{u_{k+\frac{1}{2}}^{j-1}}+\overline{u_{k+\frac{1}{2}}^{j}}{u}_{k+\frac{1}{2}}^{j-1}\right]\\
  = & \frac{2}{\tau}h\sum\limits_k\left[\left(\delta_{\hat{t}}u_{k+\frac{1}{2}}^{j+\frac{1}{2}}\overline{u_{k+\frac{1}{2}}^{j+\frac{1}{2}}}
      -\delta_{\hat{t}}\overline{u_{k+\frac{1}{2}}^{j+\frac{1}{2}}}{u}_{k+\frac{1}{2}}^{j+\frac{1}{2}}\right)
      -\left(\delta_{\hat{t}}u_{k+\frac{1}{2}}^{j-\frac{1}{2}}\overline{u_{k+\frac{1}{2}}^{j-\frac{1}{2}}}
      -\delta_{\hat{t}}\overline{u_{k+\frac{1}{2}}^{j-\frac{1}{2}}}{u}_{k+\frac{1}{2}}^{j-\frac{1}{2}}\right)\right].
 \end{array}
\end{align}
The second term
\begin{align*}
     -\frac{1}{{2}}h\sum\limits_k\delta_x^2\left(u_{k}^{j+\frac{1}{2}}+u_{k}^{j-\frac{1}{2}}\right)B\overline{u_{k}^{j}}
  =  -\frac{1}{{4}}h\sum\limits_k\delta_x^2Bu_{k}^{j}B\overline{u_{k}^{j}}
  =  \frac{1}{{4}}h\sum\limits_k\left|\delta_{\bar{x}}(Bu_{k}^{j})\right|^2,
\end{align*}
is real. The fourth term
\begin{align*}
    -\frac{i\theta}{{4}}h\sum\limits_k\delta_{2x}Bu_{k}^{j}B\overline{u_{k}^{j}}
  & =-\frac{i\theta}{{8h}}h\sum\limits_k\left(Bu_{k+1}^{j}B\overline{u_{k}^{j}}-Bu_{k-1}^{j}B\overline{u_k^{j}}\right),
\end{align*}
and the sixth term
\begin{align*}
    & \frac{\lambda}{{4}}h\sum\limits_k\left[(u_{k+\frac{1}{2}}^{j+\frac{1}{2}}+u_{k+\frac{1}{2}}^{j-\frac{1}{2}})+(
      u_{k-\frac{1}{2}}^{j+\frac{1}{2}}+u_{k-\frac{1}{2}}^{j-\frac{1}{2}})\right]
      2\left[\overline{u_k^{j+\frac{1}{2}}}+\overline{u_{k}^{j-\frac{1}{2}}}\right]\\
  = & \lambda h\sum\limits_k\left|u_{k+\frac{1}{2}}^{j+\frac{1}{2}}+u_{k+\frac{1}{2}}^{j-\frac{1}{2}}\right|^2,
\end{align*}
are real, too. The third term is
\begin{align*}
     & {i\alpha}h\sum\limits_k\left(\delta_{2t}{u_{k+\frac{1}{2}}^{j}}+\delta_{2t}u_{k-\frac{1}{2}}^{j}\right)
       \left(\overline{u_k^{j+\frac{1}{2}}}+\overline{u_k^{j-\frac{1}{2}}}\right)\\
  =  & -\frac{2i\alpha}{\tau}h\sum\limits_k\left(u_{k+\frac{1}{2}}^{j+\frac{1}{2}}-u_{k+\frac{1}{2}}^{j-\frac{1}{2}}\right)
       \left(\overline{u_{k+\frac{1}{2}}^{j+\frac{1}{2}}}+\overline{u_{k+\frac{1}{2}}^{j-\frac{1}{2}}}\right).
 \end{align*}
Abstract the imaginary part from this term, one obtains
\begin{align}\label{mass_third}
  -\frac{2\alpha i}{\tau}\left(\left\|u^{j+\frac{1}{2}}\right\|_{\frac{1}{2}}^2-\left\|u^{j-\frac{1}{2}}\right\|_{\frac{1}{2}}^2\right).
 \end{align}
 Now, we analyze the fifth term
\begin{align*}
   & {\gamma}\delta_{2t}\delta_{2x}u_k^j (\overline{u_k^{j+1}}+2\overline{u_k^{j}}+\overline{u_k^{j-1}})  \\
 = & \frac{2\gamma}{h\tau}h\sum\limits_k\left[(u_{k+\frac{1}{2}}^{j+\frac{1}{2}}-u_{k+\frac{1}{2}}^{j-\frac{1}{2}})
     -(u_{k-\frac{1}{2}}^{j+\frac{1}{2}}-u_{k-\frac{1}{2}}^{j-\frac{1}{2}})\right]
      (\overline{u_k^{j+\frac{1}{2}}}+\overline{u_k^{j-\frac{1}{2}}})\\
 = &  -\frac{2\gamma}{{h}\tau}h\sum\limits_k\left(u_{k+1}^{j+\frac{1}{2}}-u_{k+1}^{j-\frac{1}{2}}\right)
      \left[\left(\overline{u_{k+1}^{j+\frac{1}{2}}}+\overline{u_{k+1}^{j-\frac{1}{2}}}\right)
      -\left(\overline{u_{k}^{j+\frac{1}{2}}}+\overline{u_{k}^{j-\frac{1}{2}}}\right)\right]\\
 = & -\frac{2\gamma}{\tau}h\sum\limits_k\left(u_{k+1}^{j+\frac{1}{2}}-u_{k+1}^{j-\frac{1}{2}}\right)
     \left(\delta_{\hat{x}}\overline{u_{k+\frac{1}{2}}^{j+\frac{1}{2}}}+\delta_{\hat{x}}\overline{u_{k+\frac{1}{2}}^{j-\frac{1}{2}}}\right)\\
 = & -\frac{2\gamma}{\tau}h\sum\limits_k\left(u_{k+\frac{1}{2}}^{j+\frac{1}{2}}\delta_{\hat{x}}\overline{u_{k+\frac{1}{2}}^{j+\frac{1}{2}}}
     -u_{k+\frac{1}{2}}^{j-\frac{1}{2}}\delta_{\hat{x}}\overline{u_{k+\frac{1}{2}}^{j-\frac{1}{2}}}\right)\\
   & -\frac{2\gamma}{\tau}h\sum\limits_k\left(u_{k+\frac{1}{2}}^{j+\frac{1}{2}}\delta_{\hat{x}}\overline{u_{k+\frac{1}{2}}^{j-\frac{1}{2}}}
     -u_{k+\frac{1}{2}}^{j-\frac{1}{2}}\delta_{\hat{x}}\overline{u_{k+\frac{1}{2}}^{j+\frac{1}{2}}}\right).
\end{align*}
Based on the second conclusion in Lemma \ref{Green}, the  first part of the above equality is
purely imaginary, and the second part is real. Therefore, the
imaginary part of the term is
\begin{align}\label{mass_fifth}
-\frac{2\gamma}{\tau}h\sum\limits_k\left(u_{k+\frac{1}{2}}^{j+\frac{1}{2}}\delta_{\hat{x}}\overline{u_{k+\frac{1}{2}}^{j+\frac{1}{2}}}
     -u_{k+\frac{1}{2}}^{j-\frac{1}{2}}\delta_{\hat{x}}\overline{u_{k+\frac{1}{2}}^{j-\frac{1}{2}}}\right).
\end{align}
Finally, the last term is
\begin{align*}
   & \frac{2\beta}{4}h\sum\limits_k\left(\left|u_{k+\frac{1}{2}}^{j+\frac{1}{2}}\right|^2u_{k+\frac{1}{2}}^{j+\frac{1}{2}}
      +\left|u_{k+\frac{1}{2}}^{j-\frac{1}{2}}\right|^2u_{k+\frac{1}{2}}^{j-\frac{1}{2}}
      +\left|u_{k-\frac{1}{2}}^{j+\frac{1}{2}}\right|^2u_{k-\frac{1}{2}}^{j+\frac{1}{2}}
      +\left|u_{k-\frac{1}{2}}^{j-\frac{1}{2}}\right|^2u_{k-\frac{1}{2}}^{j-\frac{1}{2}}\right)
      \left(\overline{u_{k}^{j+\frac{1}{2}}}+\overline{u_{k}^{j-\frac{1}{2}}}\right)\\
 =  &\frac{\beta}{2} h\sum\limits_k\left(\left|u_{k+\frac{1}{2}}^{j+\frac{1}{2}}\right|^2u_{k+\frac{1}{2}}^{j+\frac{1}{2}}
      +\left|u_{k+\frac{1}{2}}^{j-\frac{1}{2}}\right|^2u_{k+\frac{1}{2}}^{j-\frac{1}{2}}\right)
      \left(\overline{u_{k+\frac{1}{2}}^{j+\frac{1}{2}}}+\overline{u_{k+\frac{1}{2}}^{j-\frac{1}{2}}}\right).
\end{align*}
The imaginary part this term is
\begin{align}\label{mass_last}
 \begin{array}{ll}
   & \frac{\beta}{2}\mathcal{I}h\sum\limits_k\left(\left|u_{k+\frac{1}{2}}^{j+\frac{1}{2}}\right|u_{k+\frac{1}{2}}^{j+\frac{1}{2}}
     \overline{u_{k+\frac{1}{2}}^{j-\frac{1}{2}}}+
     \left|u_{k+\frac{1}{2}}^{j-\frac{1}{2}}\right|u_{k+\frac{1}{2}}^{j-\frac{1}{2}}\overline{u_{k+\frac{1}{2}}^{j+\frac{1}{2}}}\right)\\
 = & \frac{\beta}{4}h\sum\limits_k\left[\left|u_{k+\frac{1}{2}}^{j+\frac{1}{2}}\right|^2-\left|u_{k+\frac{1}{2}}^{j-\frac{1}{2}}\right|^2\right]
     \left[u_{k+\frac{1}{2}}^{j+\frac{1}{2}}\overline{u_{k+\frac{1}{2}}^{j-\frac{1}{2}}}-u_{k+\frac{1}{2}}^{j-\frac{1}{2}}\overline{u_{k+\frac{1}{2}}^{j+\frac{1}{2}}}\right]\\
 = & \frac{\beta}{4}h\sum\limits_k\left[\left|u_{k+\frac{1}{2}}^{j+\frac{1}{2}}\right|^2-\left|u_{k+\frac{1}{2}}^{j-\frac{1}{2}}\right|^2\right]
     \left(u_{k+\frac{1}{2}}^{j+\frac{1}{2}}-u_{k+\frac{1}{2}}^{j-\frac{1}{2}}\right)
     \left(\overline{u_{k+\frac{1}{2}}^{j+\frac{1}{2}}}+\overline{u_{k+\frac{1}{2}}^{j-\frac{1}{2}}}\right)\\
   & -\frac{\beta}{4}h\sum\limits_k\left[\left|u_{k+\frac{1}{2}}^{j+\frac{1}{2}}\right|^2-\left|u_{k+\frac{1}{2}}^{j-\frac{1}{2}}\right|^2\right]^2.
 \end{array}
\end{align}
Summing over the equalities (\ref{mass_first}), (\ref{mass_third}), (\ref{mass_fifth}) and (\ref{mass_last}),
 one obtains the conclusion what one wishes.
\end{pro}

\section{Numerical examples}\label{example}
In this section, we investigate the theoretical analysis by  a
series of numerical experiments, including conservation properties
and the accuracy of the schemes.

\begin{figure}[h]
   \begin{center}
       \epsfig{file=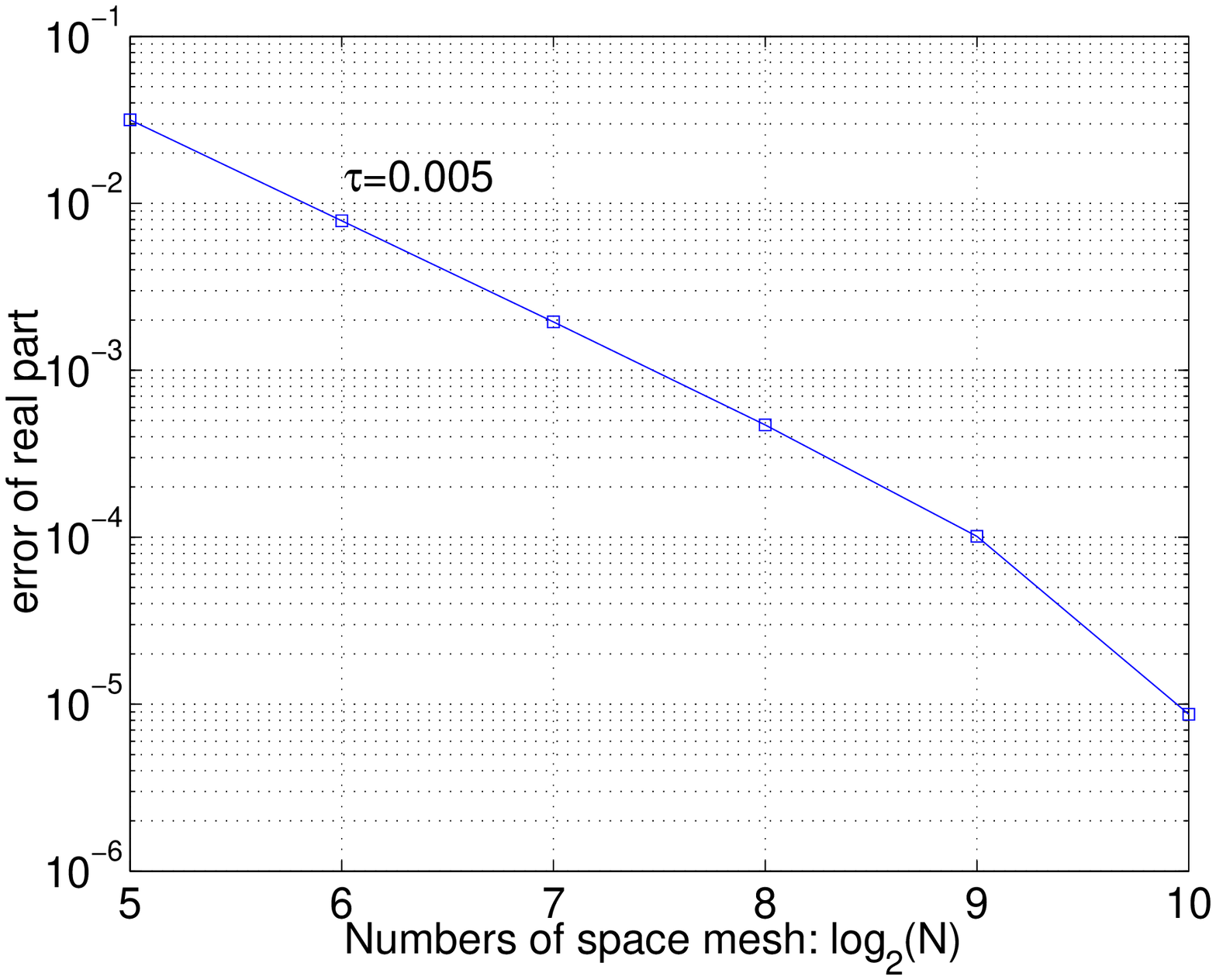,height=4.0cm,width=6.5cm}, ~~
       \epsfig{file=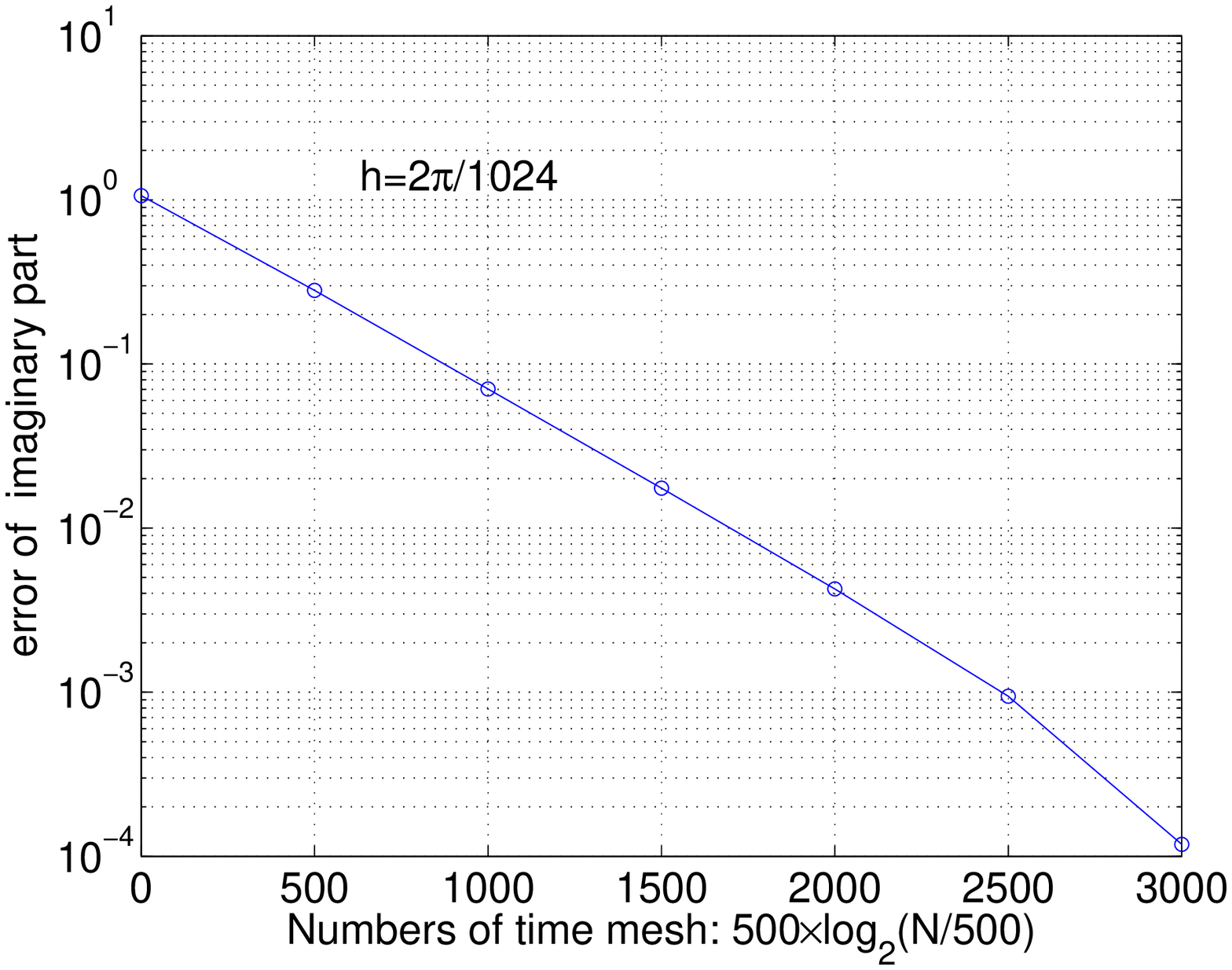,height=4.0cm,width=6.5cm}\\
       \caption{\label{space_time_order}\small Numerical error vs. mesh numbers: left for space, right for time.}
   \end{center}
\end{figure}

\begin{figure}[h]
   \begin{center}
       \epsfig{file=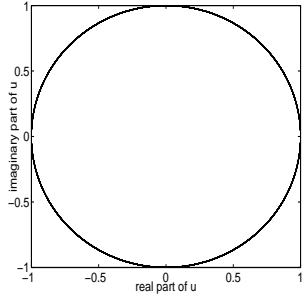,height=4.0cm,width=4.0cm}\\
       \caption{\label{phase}\small Phasic profiles of numerical solution at all time steps.}
   \end{center}
\end{figure}

\begin{figure}[h]
   \begin{center}
       \epsfig{file=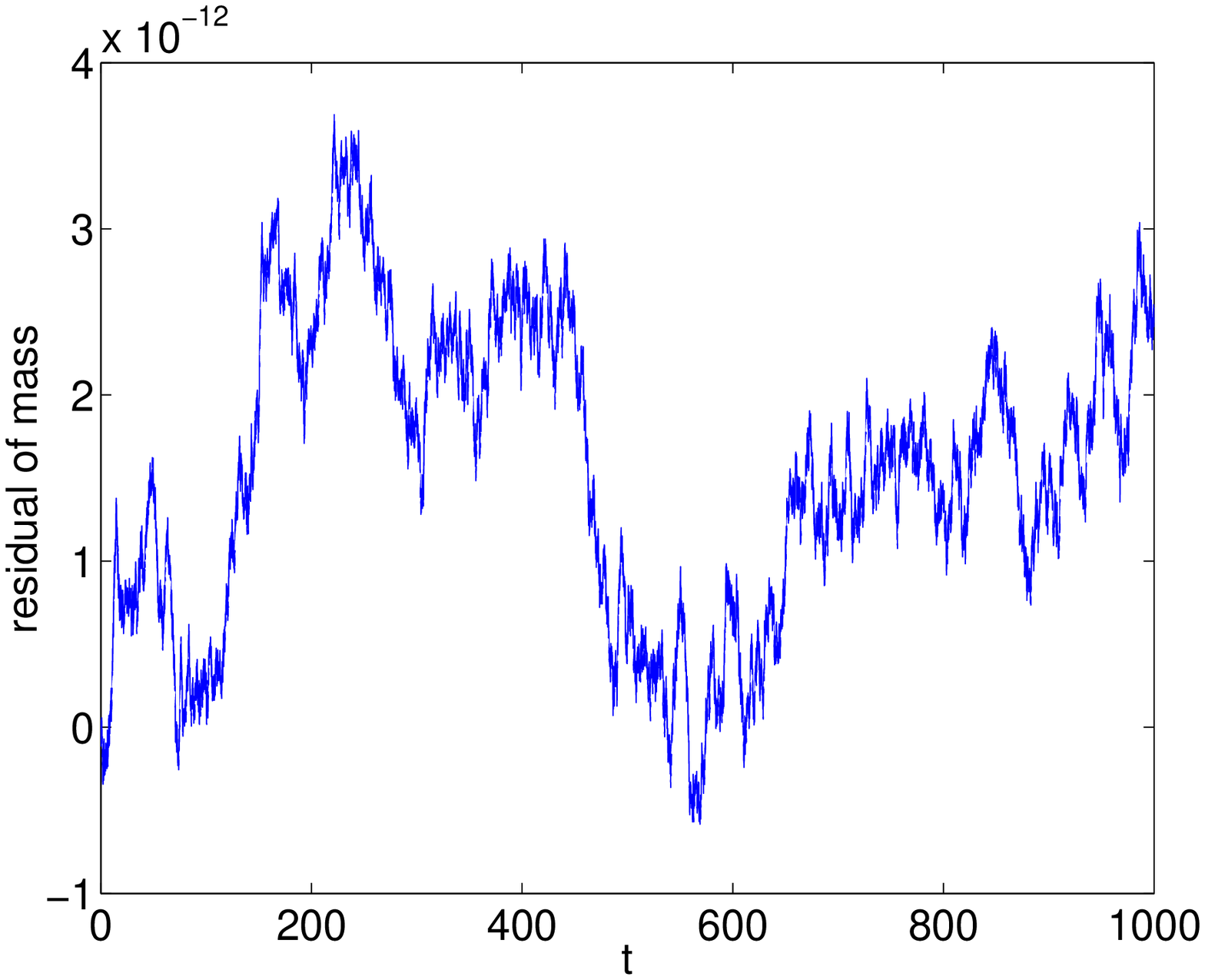,height=4.0cm,width=6.5cm},
       \epsfig{file=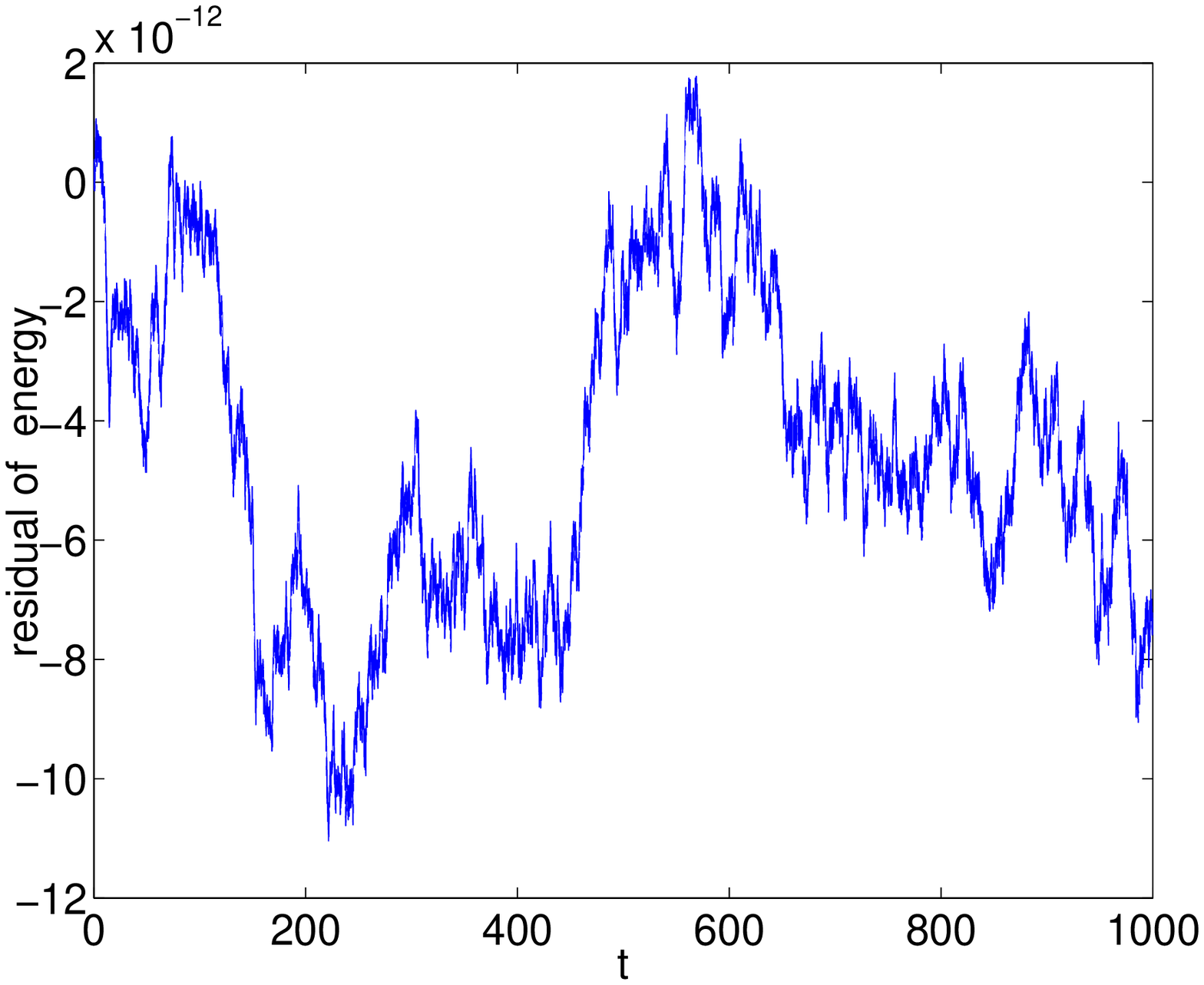,height=4.0cm,width=6.5cm}\\
       \caption{\label{linear_energy_mass}\small Residuals of mass (left) and energy (right).}
   \end{center}
\end{figure}

\begin{example}\em Firstly, we consider a linear problem, i.e., $\beta=0$,\end{example}
\begin{align}\label{linear}
  \left\{
   \begin{array}{l}
      u_{tt}-u_{xx}+u_{tx}+i(u_t+u_x)+3u=0,~(x,t)\in{[0, 2\pi]\times{(0, 50]}},\\
      u(x+2\pi,t)=u(x,t),
      u(x,0)=\exp{(ix)}, u_t(x,0)=-3i\exp{(ix)}.
   \end{array}\right.
\end{align}
In this case, the energy and mass are constants from  Theorems \ref{discrete_energy} and \ref{discrete_mass}.
The exact solution of the determined problem is
\begin{align}\label{linear_exact}
  u(x,t)=\exp(i(x-3t)),
\end{align}
which is a plane wave propagating to the right with velocity $v=3$. The amplitude of the wave is equal to $1$.
We simulate the problem by the MI (\ref{Preissman}) under diverse mesh divisions. The left of Fig.
\ref{space_time_order} shows the maximum error for the real part of numerical solution against the space mesh numbers under
$\tau=0.005$ at $t=50$, and the right one presents the maximum error for the imaginary part of numerical solution against the time mesh numbers under
$h=\frac{2\pi}{1024}$ at $t=50$. The error is metered $e=u_k^j-U_k^j$. The figures imply that the numerical solution of the MI (\ref{Preissman})
 converges to the exact solution almost with the same rate $2$ both in time and space.
Fig. \ref{phase} presents the phasic profiles of the numerical solution at every time step. It
is suggested that the curves are overlapped and always keep a unit circle. This exactly agrees with the exact solution (\ref{linear_exact}).
Fig. \ref{linear_energy_mass} plots the residuals of energy (left) and mass (right) with $\tau=0.01, h=\frac{2\pi}{64}$
against time $t\in[0,1000]$. Judging from the plots, the MI preserves both the energy and mass indeed in case of $\beta=0$.

\begin{figure}[h]
   \begin{center}
       \epsfig{file=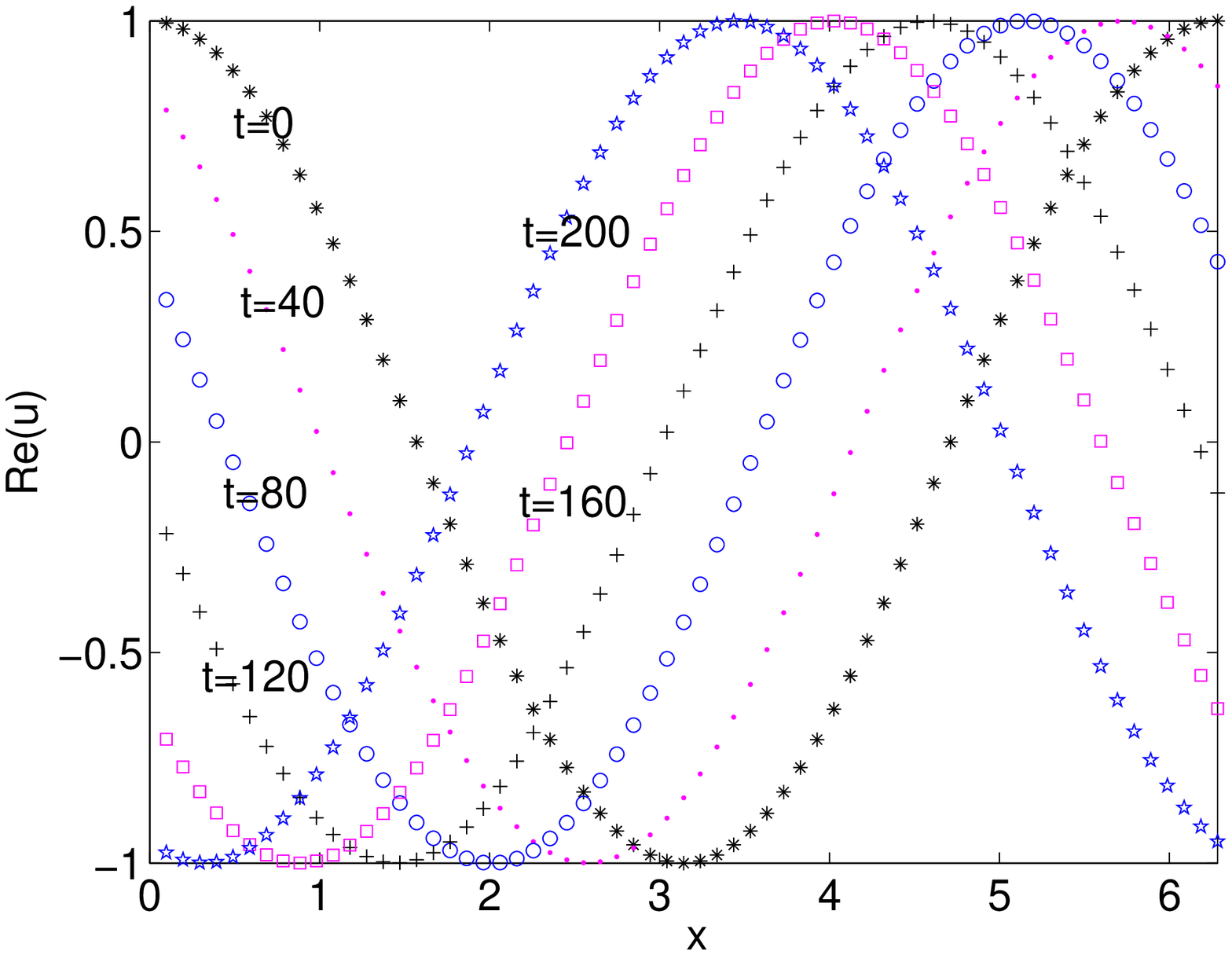,height=4.0cm,width=6.5cm},
       \epsfig{file=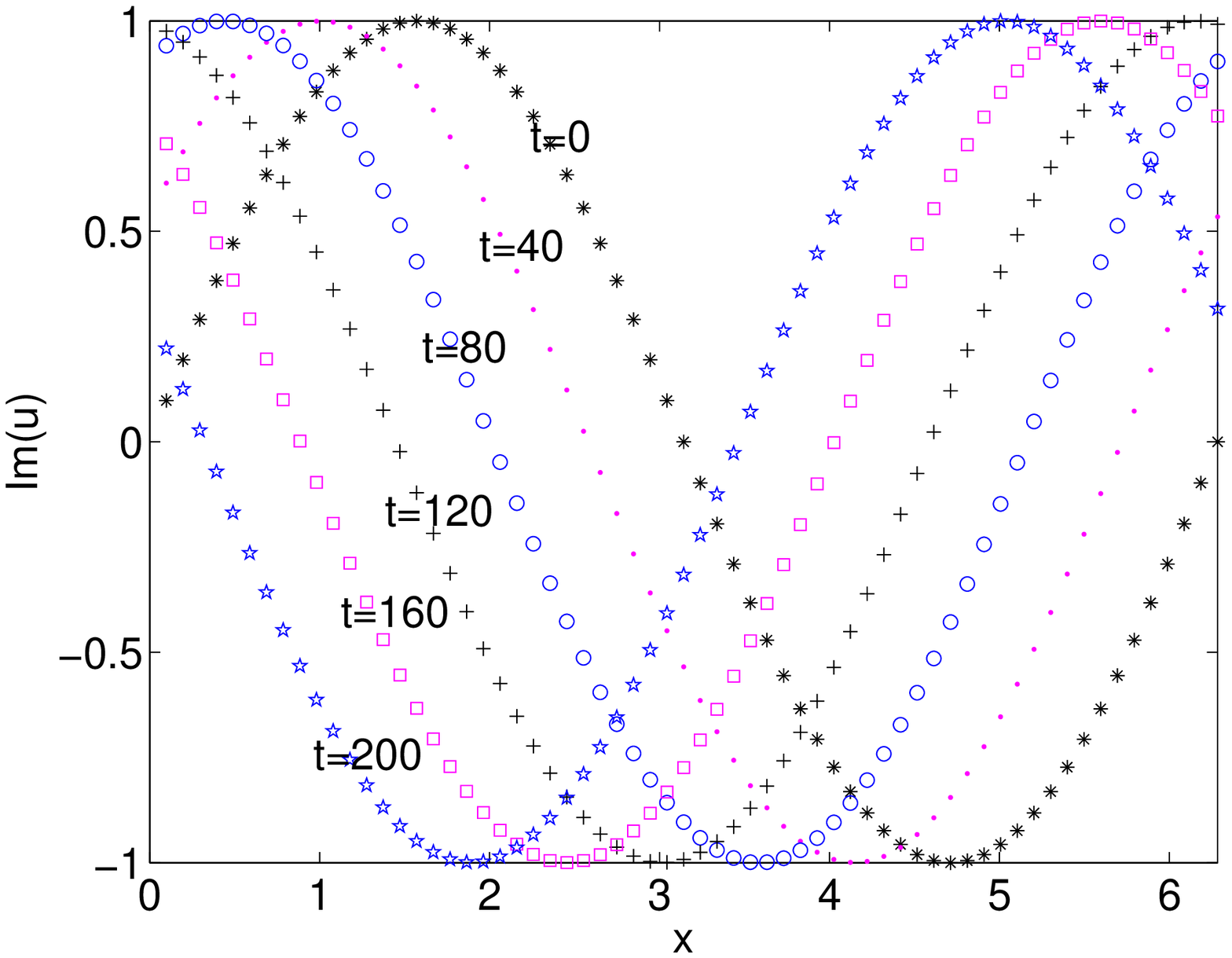,height=4.0cm,width=6.5cm}\\
       \caption{\label{nonlinear_periodic_u}\small The real (left) and imaginary (right)
       parts of the numerical solution at different times.}
   \end{center}
\end{figure}

\begin{figure}[h]
   \begin{center}
       \epsfig{file=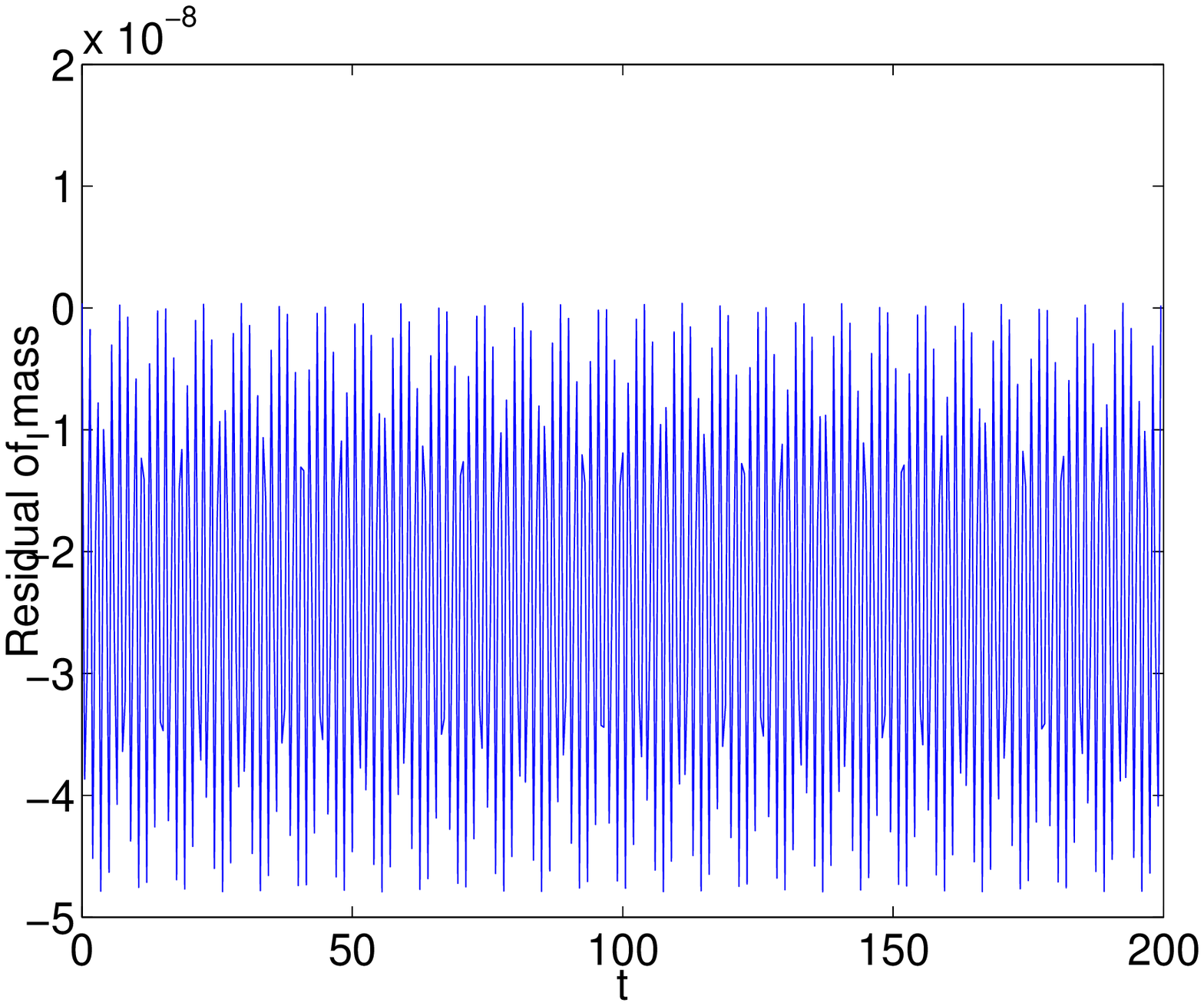,height=4.0cm,width=6.5cm},
       \epsfig{file=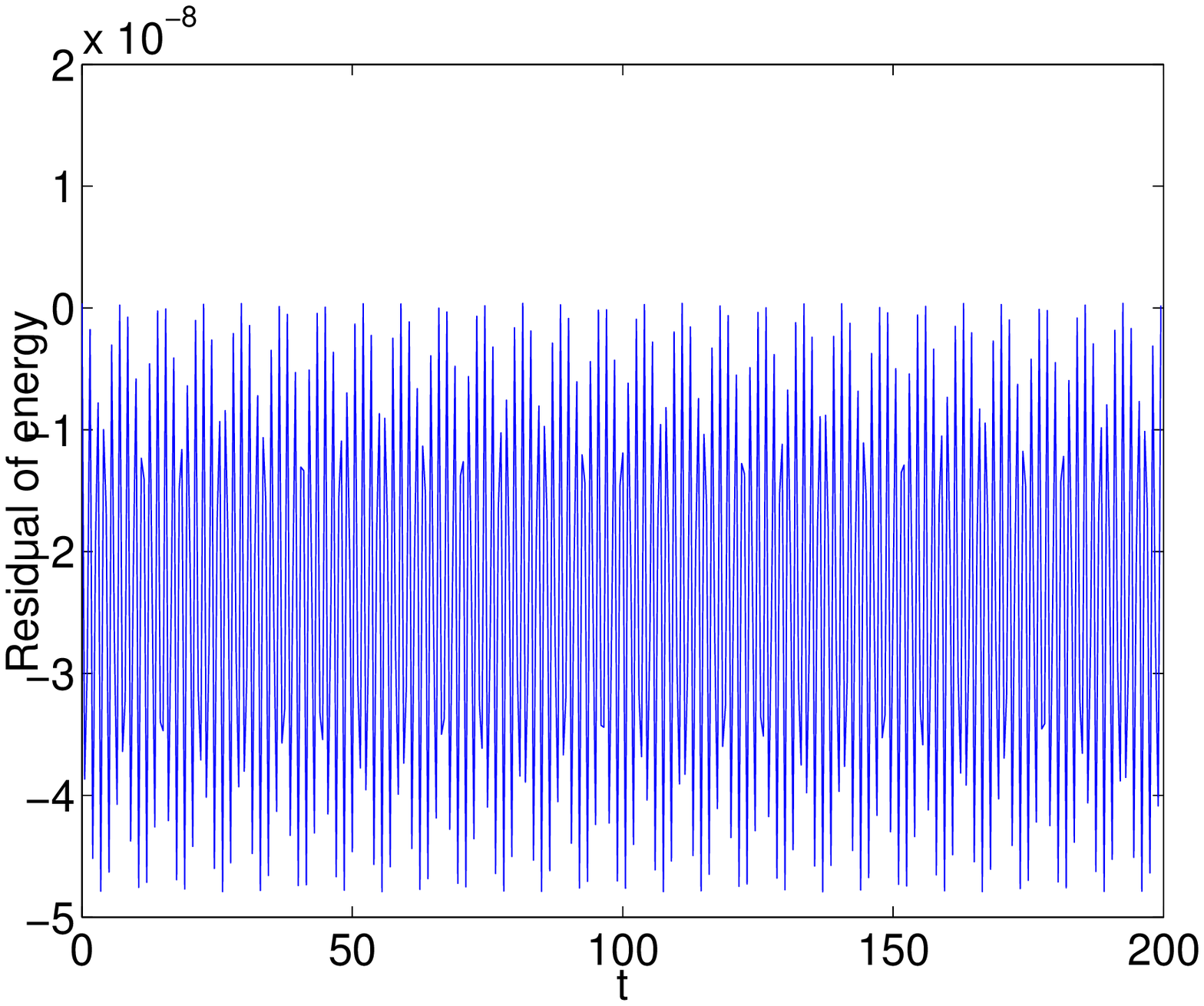,height=4.0cm,width=6.5cm}\\
       \caption{\label{nonlinear_periodic_redidual}\small The residuals of mass (left) and energy (right).}
   \end{center}
\end{figure}

\begin{example}\em We simulate the following problem\end{example}
\begin{align}\label{Schr_wave}
  \left\{
 \begin{array}{l}
   u_{tt}-u_{xx}+u_{tx}+i(u_t+u_x)+{u}+2{\left|u\right|^2}u=0,~(x,t)\in{[0, 2\pi]\times{(0, 200]}},\\
   u(x+2\pi,t)=u(x,t),~~
   u(x,0)=\exp(ix),~~u_t(x,0)=i\exp(ix),
 \end{array}\right.
\end{align}
by MI (\ref{Preissman}) until $T=200$.
The problem admits the following exact solution
\begin{align}
  u(x,t)=\exp(i(x+t)).
\end{align}
The spatial-temporal domain is divided by $\tau=0.01, h=\frac{\pi}{100}$. The real and imaginary part of the numerical solution at different
time are profiled in Fig. \ref{nonlinear_periodic_u}, and the residuals of mass and energy against time are presented in Fig. \ref{nonlinear_periodic_redidual}. From the figures, we can find that the curves of the real part always
follow cosine evolution, and the imaginary part are sine. The residuals of mass and energy take on periodic evolution.
It is very interesting that the residual plots of the mass and energy are very like.

\begin{example}\label{periodic} \em Next, we test the periodic initial valued problem\end{example}
\begin{align}\label{plane_3}
  \left\{
 \begin{array}{l}
   u_{tt}-u_{xx}+iu_t+2{\left|u\right|^2}u=0,~(x,t)\in{[0, 2\pi]\times{(0, 100]}},\\
   u(x+2\pi,t)=u(x,t),\\
   u(x,0)=\sqrt{3}\exp(6ix),~~u_t(x,0)=-7\sqrt{3}i\exp(6ix),
 \end{array}\right.
\end{align}
We simulate the problem by the MI (\ref{Preissman}) and Wang's energy preserving scheme
\begin{align}\label{scheme_ZW}
  \delta_t^2u_k^j-\frac{1}{2}(\delta_x^2u_{k}^{j+1}+\delta_x^2u_k^{j-1})+i\alpha\delta_{\hat{t}}u_k^j+
  (\left|u_k^{j+1}\right|^2+\left|u_k^{j-1}\right|^2)\frac{u_k^{j+1}+u_k^{j-1}}{2}=0.
\end{align}
This scheme preserves the energy (\ref{energy}) exactly \cite{WTC:AMC:06}, that is,
\begin{align}
  \mathcal{E}^{n}=\left\|(u_k^j)_t\right\|^2+\frac{1}{2}\left(\left\|(u_k^{j+1})_x\right\|^2+\left\|(u_k^j)_x\right\|^2\right)+
                  h\sum\limits_k\left|u_k^j\right|^4=\mathcal{E}^0.
\end{align}

\begin{figure}[h]
   \begin{center}
       \epsfig{file=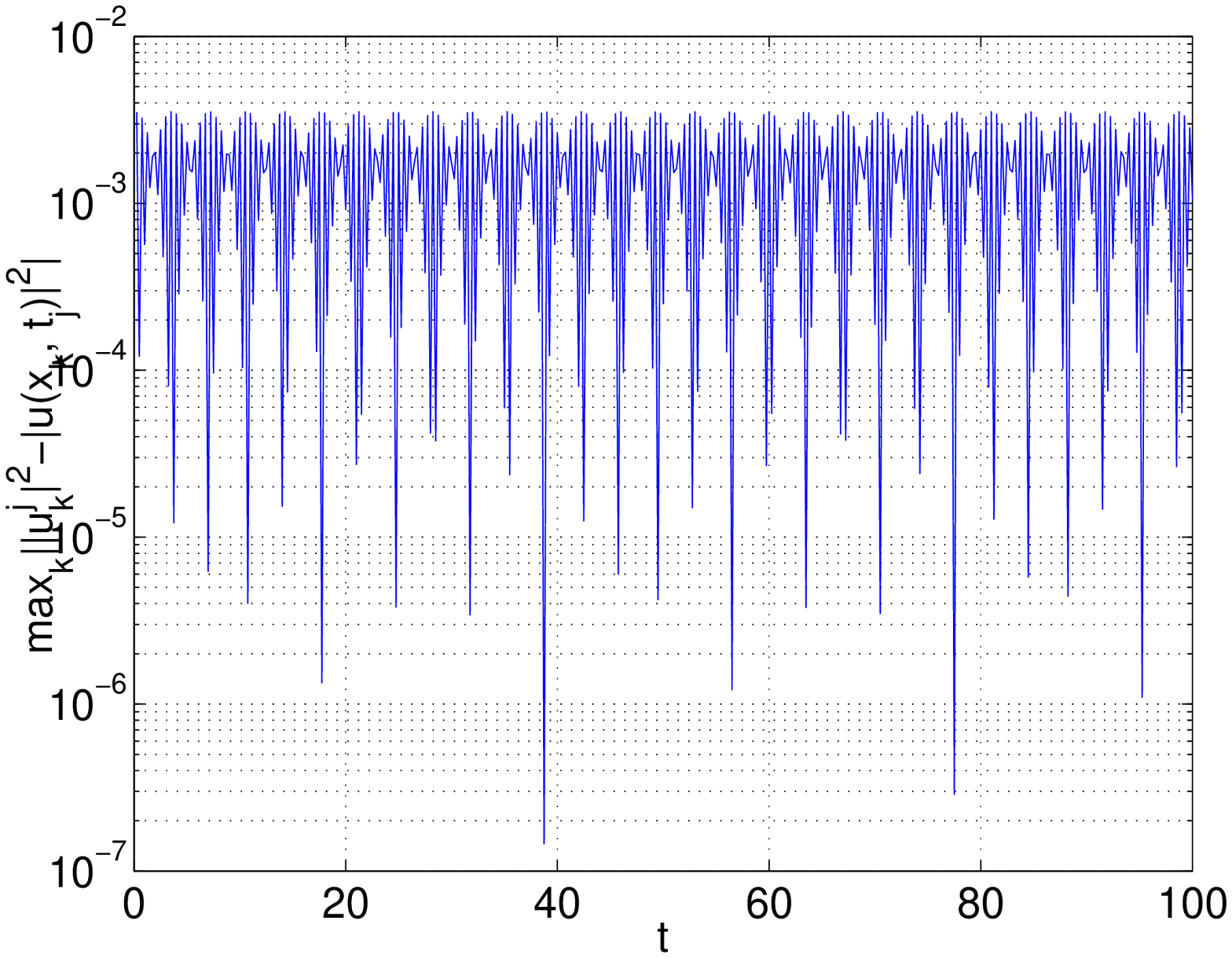,height=5cm,width=6.5cm},
       \epsfig{file=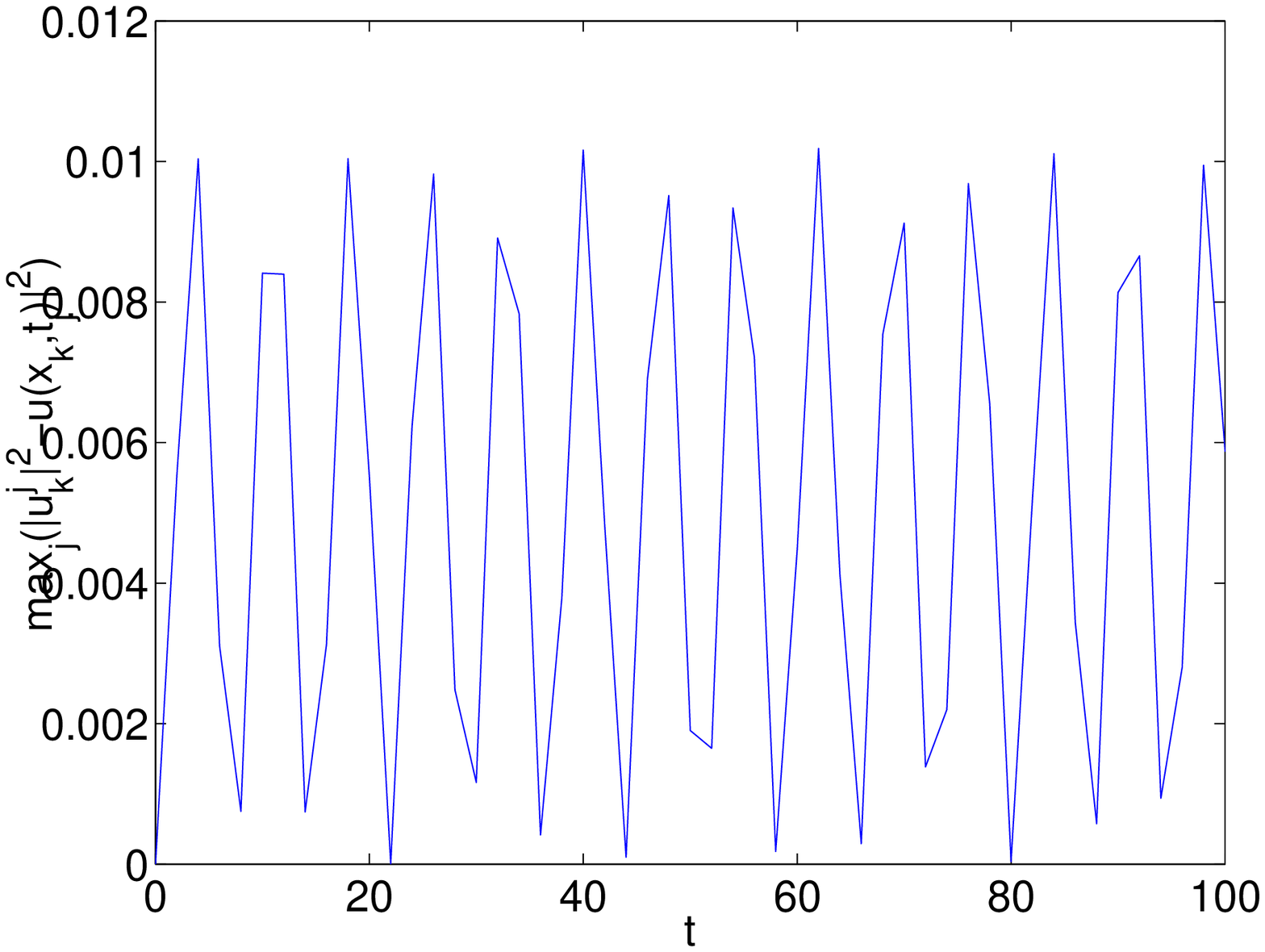,height=5cm,width=6.5cm}\\
       \caption{\label{compare_u_infty}\small Comparison $e_\infty$ by schemes
       (\ref{Preissman}) and (\ref{scheme_ZW}). Left for scheme (\ref{Preissman}), right for (\ref{scheme_ZW}).}
   \end{center}
\end{figure}

\begin{figure}
   \begin{center}
       \epsfig{file=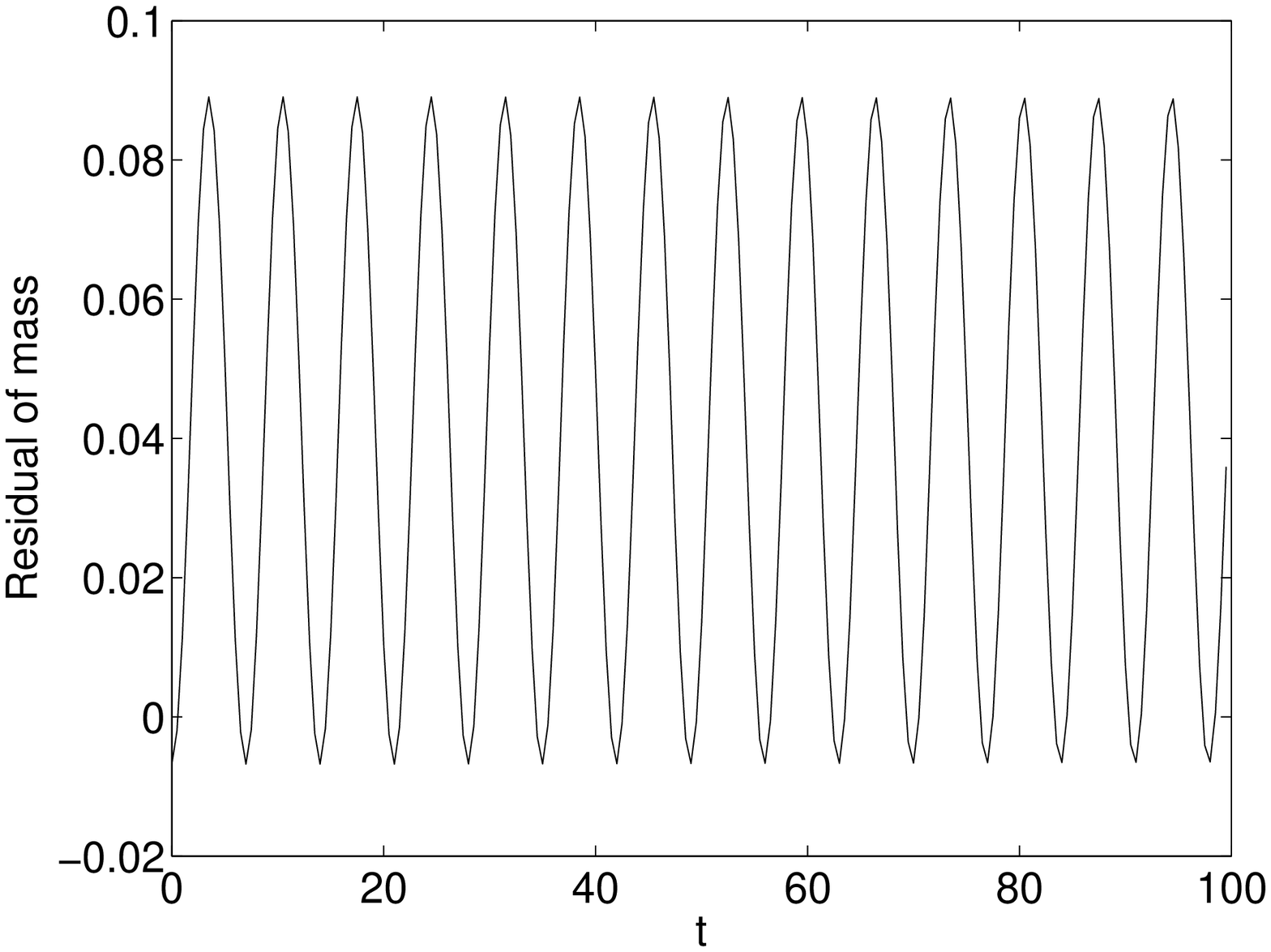,height=5cm,width=6.5cm}
       \epsfig{file=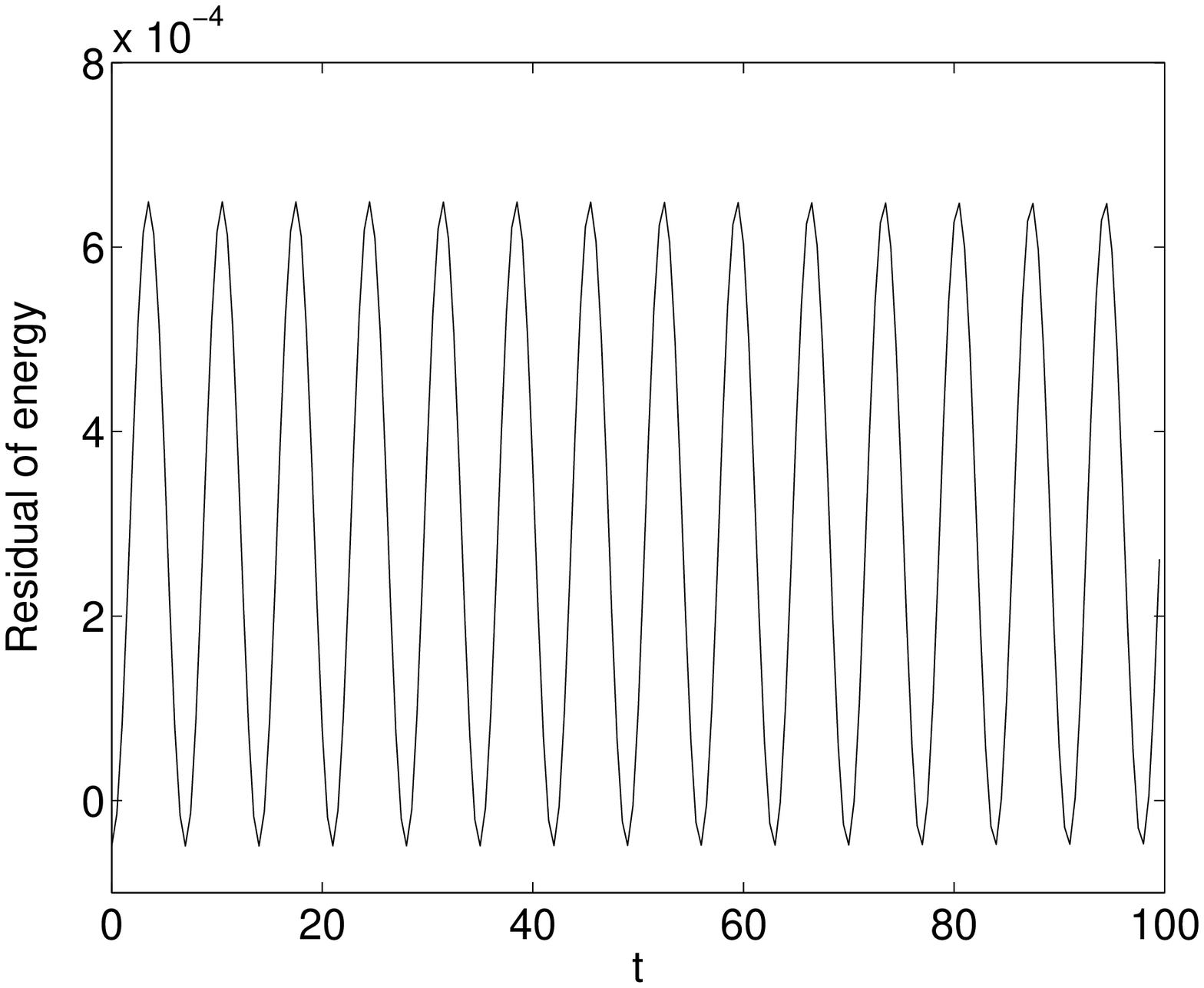,height=5cm,width=6.5cm}  \\
       \epsfig{file=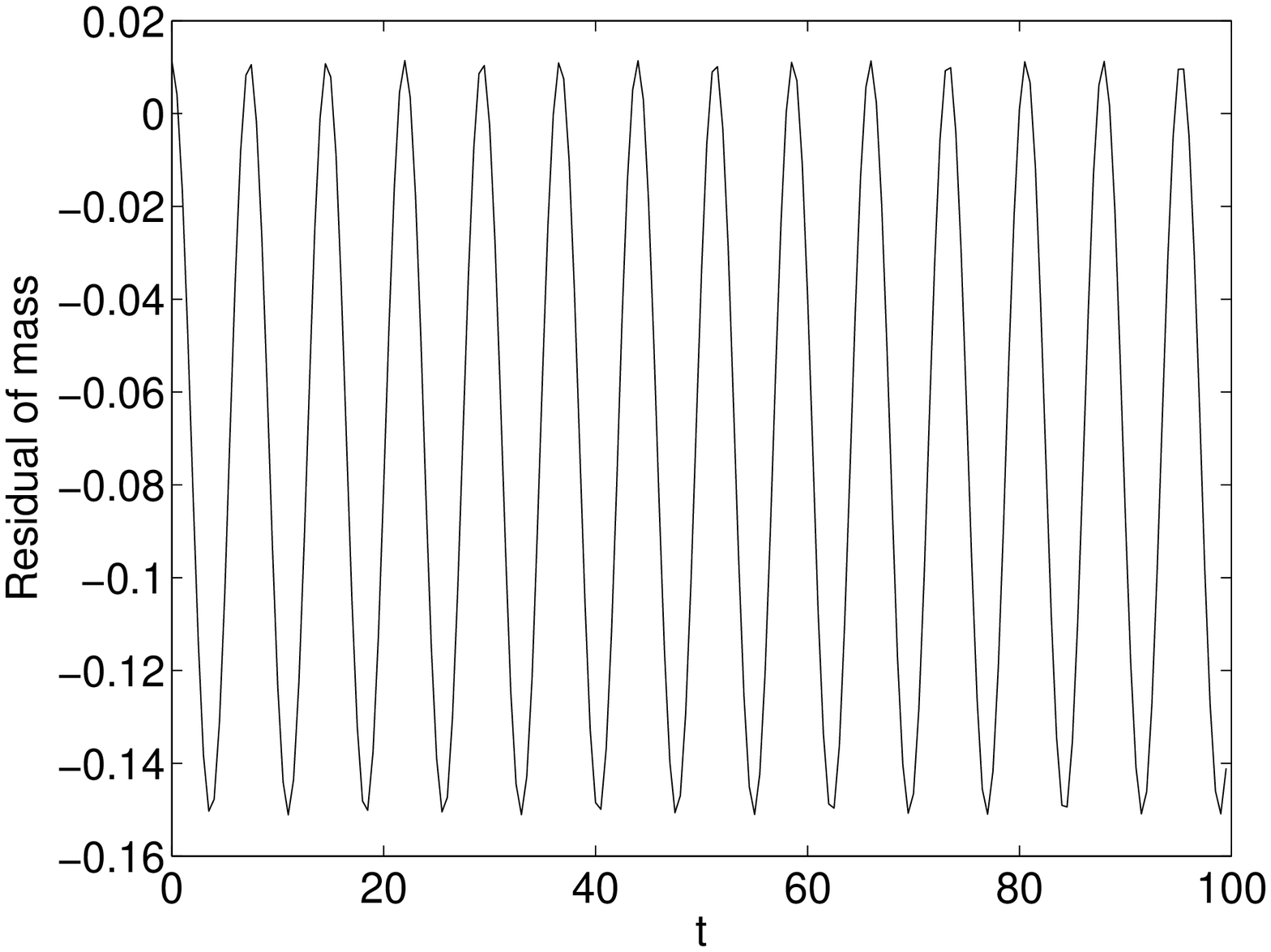,height=5cm,width=6.5cm}
       \epsfig{file=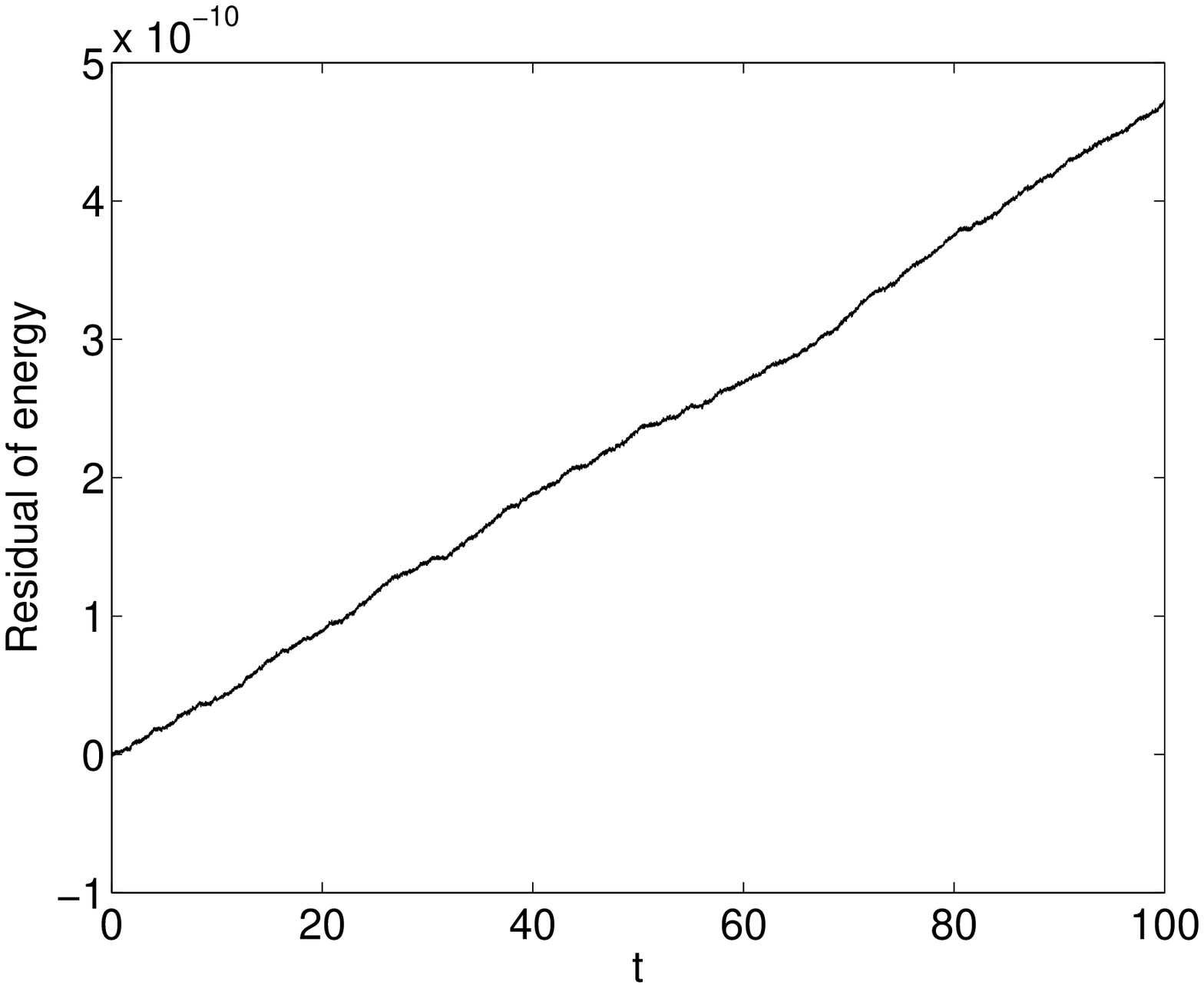,height=5cm,width=6.5cm}  \\
       \caption{\label{compare_conservative}\small Conservative comparison. Upper for scheme (\ref{Preissman}), lower for scheme (\ref{scheme_ZW}). Left for mass, right for energy.}
   \end{center}
\end{figure}

We take $T=100$, $h=\frac{2\pi}{200}, \tau=0.01$ in the example. The error of numerical solution which is measured
$e_\infty=\max\limits_k\left||u_k^j|^2-|u(x_k, t_j)|^2\right|$ by schemes (\ref{Preissman}) and (\ref{scheme_ZW}) is
presented in Fig. \ref{compare_u_infty}, and the comparison of the conservative properties is plotted in Fig. \ref{compare_conservative}.

From the figures, we can observe the following phenomena: i) The numerical solution of the multisymplectic scheme (\ref{Preissman}) is
more accurate than that of Wang's scheme (\ref{scheme_ZW}) in the test. ii) Wang's scheme (\ref{scheme_ZW}) preserves the
energy up to $10^{-10}$ scale, but not mass. The residuals of energy and mass for multisymplectic scheme (\ref{Preissman})
are fluctuated periodically, furthermore, they are very small relatively to their exact values, up to $10^{-7}$ and $10^{-4}$ scale, respectively.

\begin{example}
 \em In the example, we consider the following problem\cite{WSS:JCAM:11}
\end{example}
\begin{align}\label{plane_3}
  \left\{
 \begin{array}{l}
   u_{tt}-u_{xx}+iu_t+2{\left|u\right|^2}u=0,~(x,t)\in{[-50, 50]\times{(0,500]}},\\
   u(x,0)=Asech(Kx), ~~u_t(x,0)=i\nu Asech(Kx).
 \end{array}\right.
\end{align}
The problem has the solution in  the form
\begin{align}
  u(x,t)=Asech(Kx)\exp(i\nu t),
\end{align}
where $A=|K|, \nu=\frac{1}{2}(-1\pm \sqrt{1-4K^2})$.

We take $K=\frac{1}{4}, \nu=-\frac{1}{2}-\frac{\sqrt{3}}{4}$ in the test.

\begin{figure}[h]
   \begin{center}
       \epsfig{file=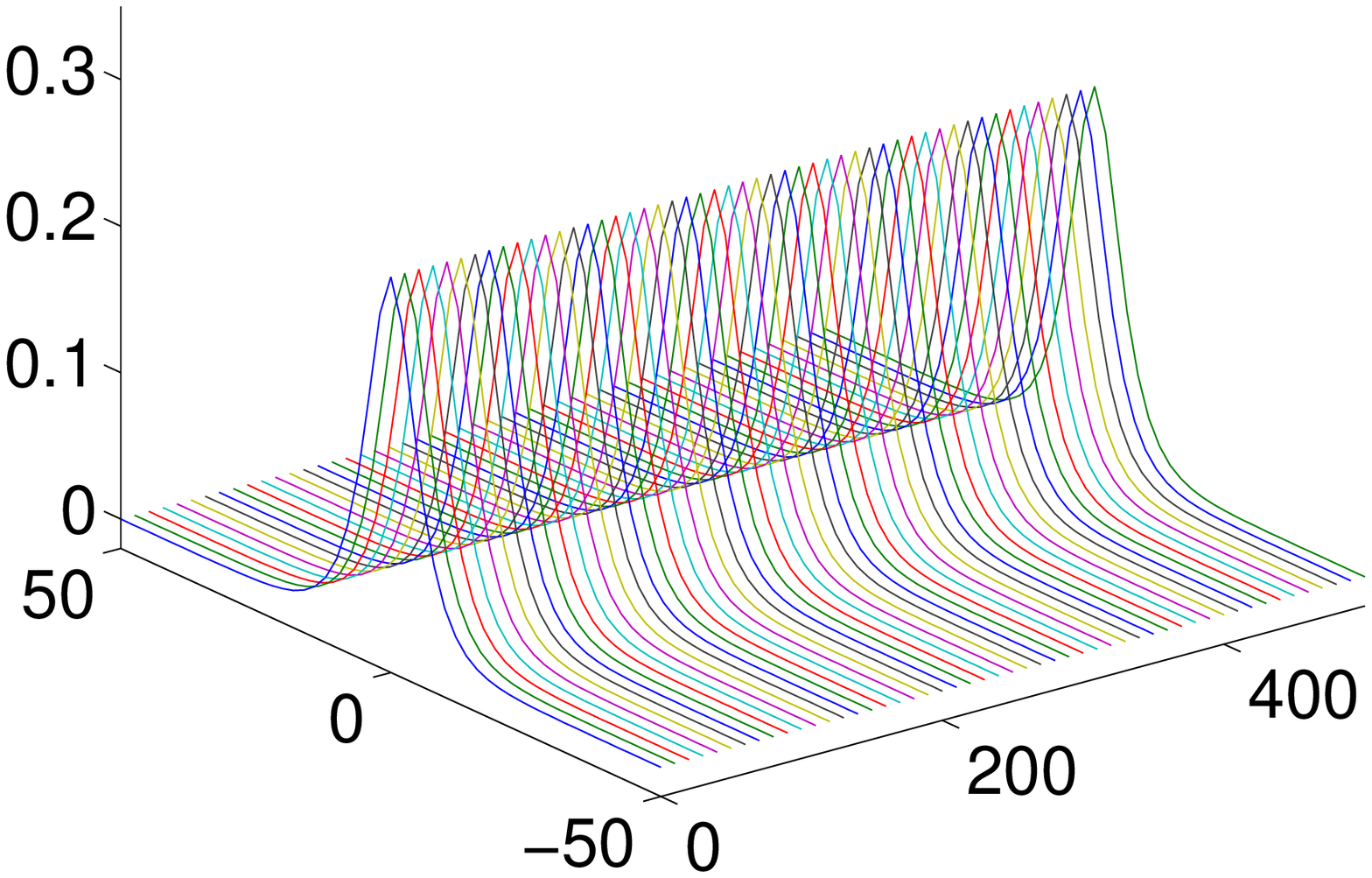,height=4.0cm,width=6.5cm},
       \epsfig{file=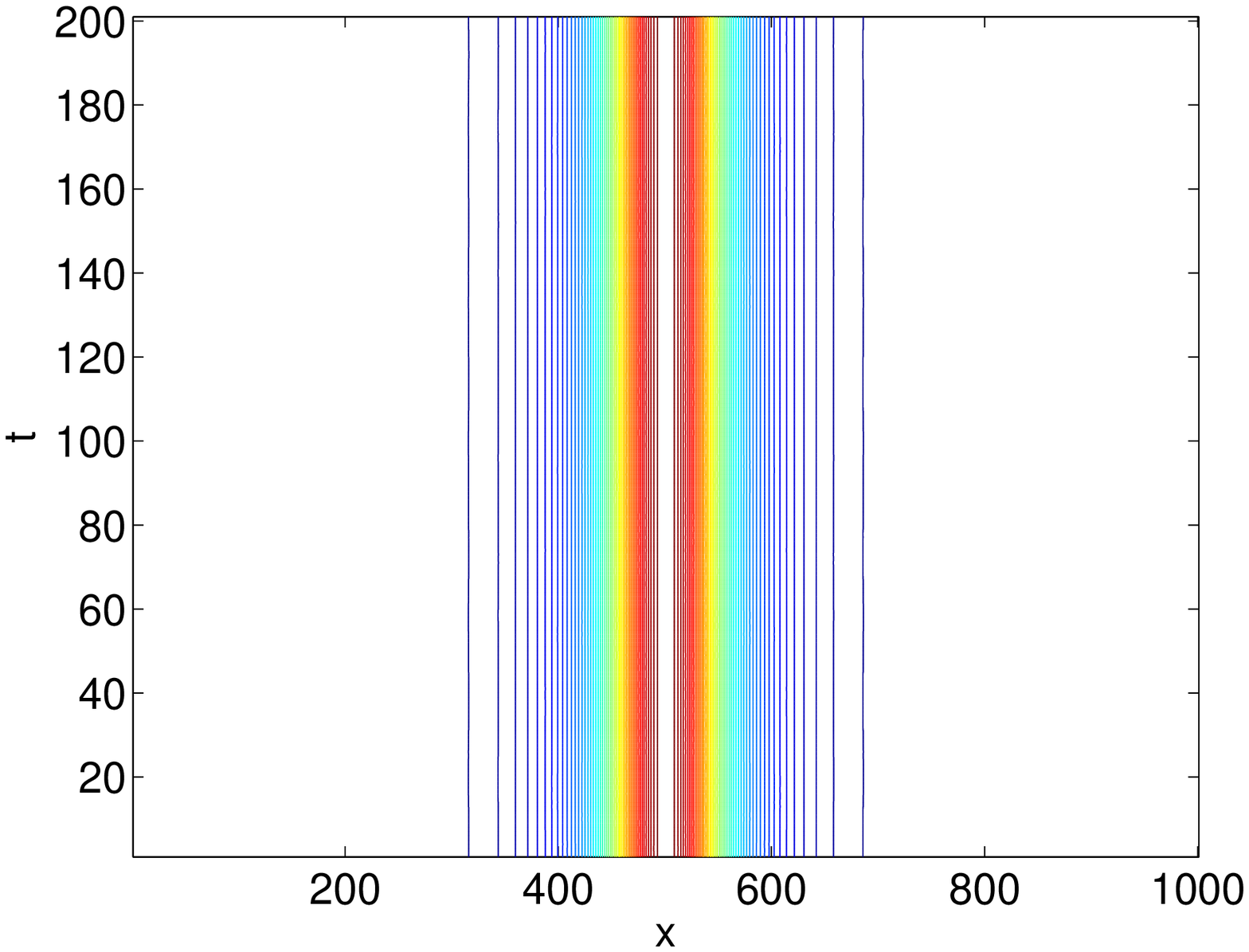,height=4.0cm,width=6.5cm}\\
       \caption{\label{abs_u}\small The solitary wave shape (left) and contours (right) of $|u|$ by scheme (\ref{Preissman}).}
   \end{center}
\end{figure}

\begin{figure}[h]
   \begin{center}
       \epsfig{file=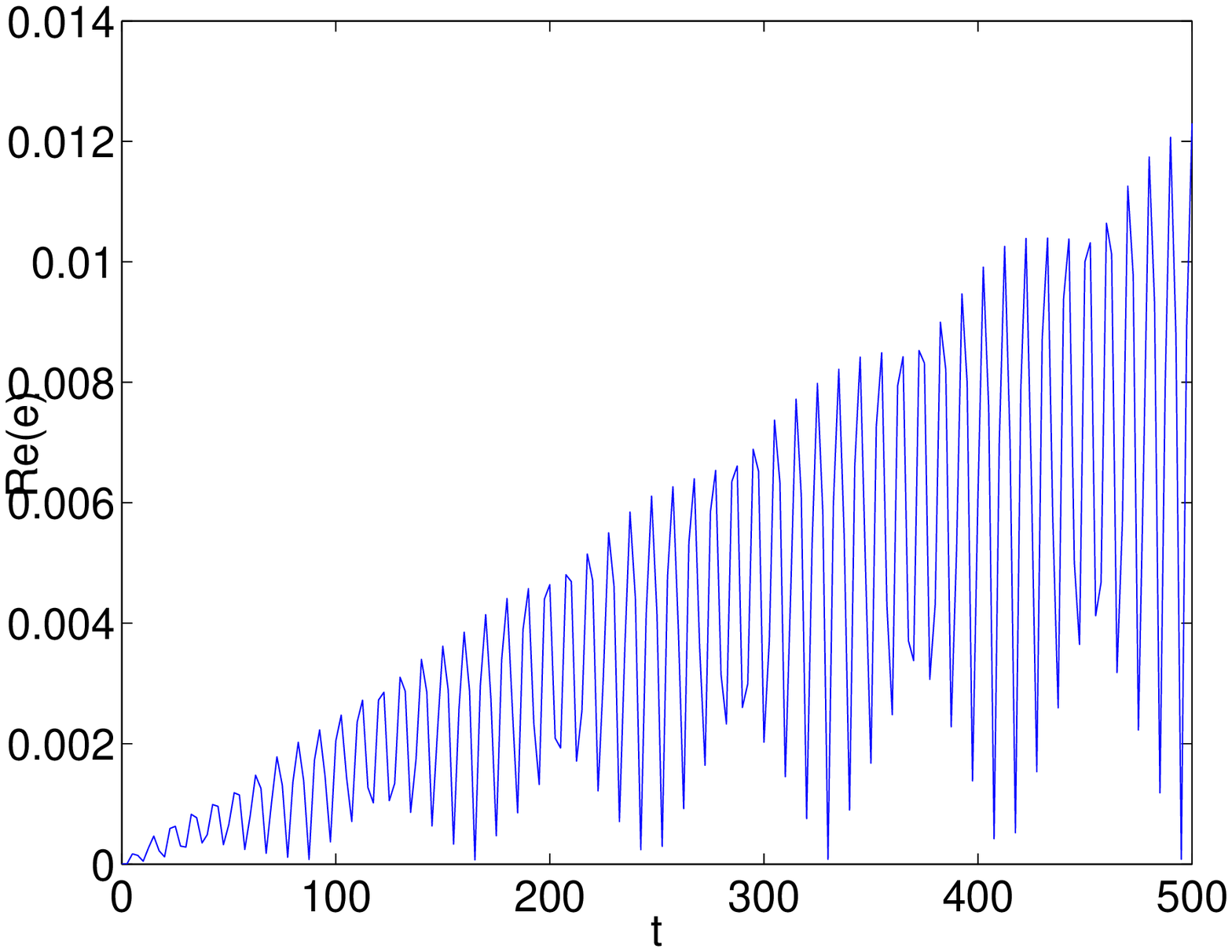,height=4.0cm,width=6.5cm},
       \epsfig{file=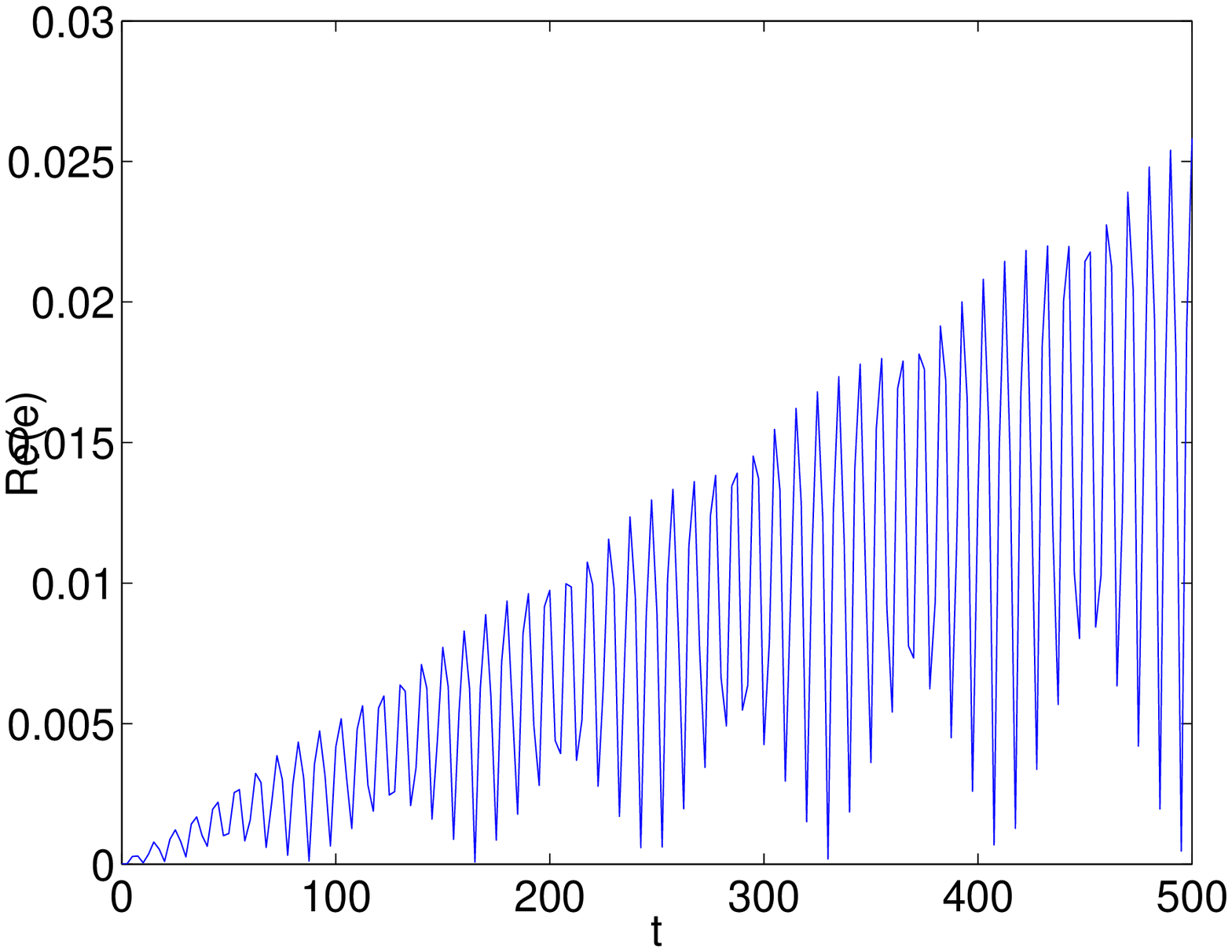,height=4.0cm,width=6.5cm}\\
       \caption{\label{compare_hom_error}\small Comparison $Re(e)=\max\limits_k\left||u_k^j|-|u(x_k, t_j)|\right|$ by schemes
       (\ref{Preissman}) and (\ref{scheme_ZW}). Left for scheme (\ref{Preissman}), right for (\ref{scheme_ZW}).}
   \end{center}
\end{figure}

\begin{figure}[h]
   \begin{center}
       \epsfig{file=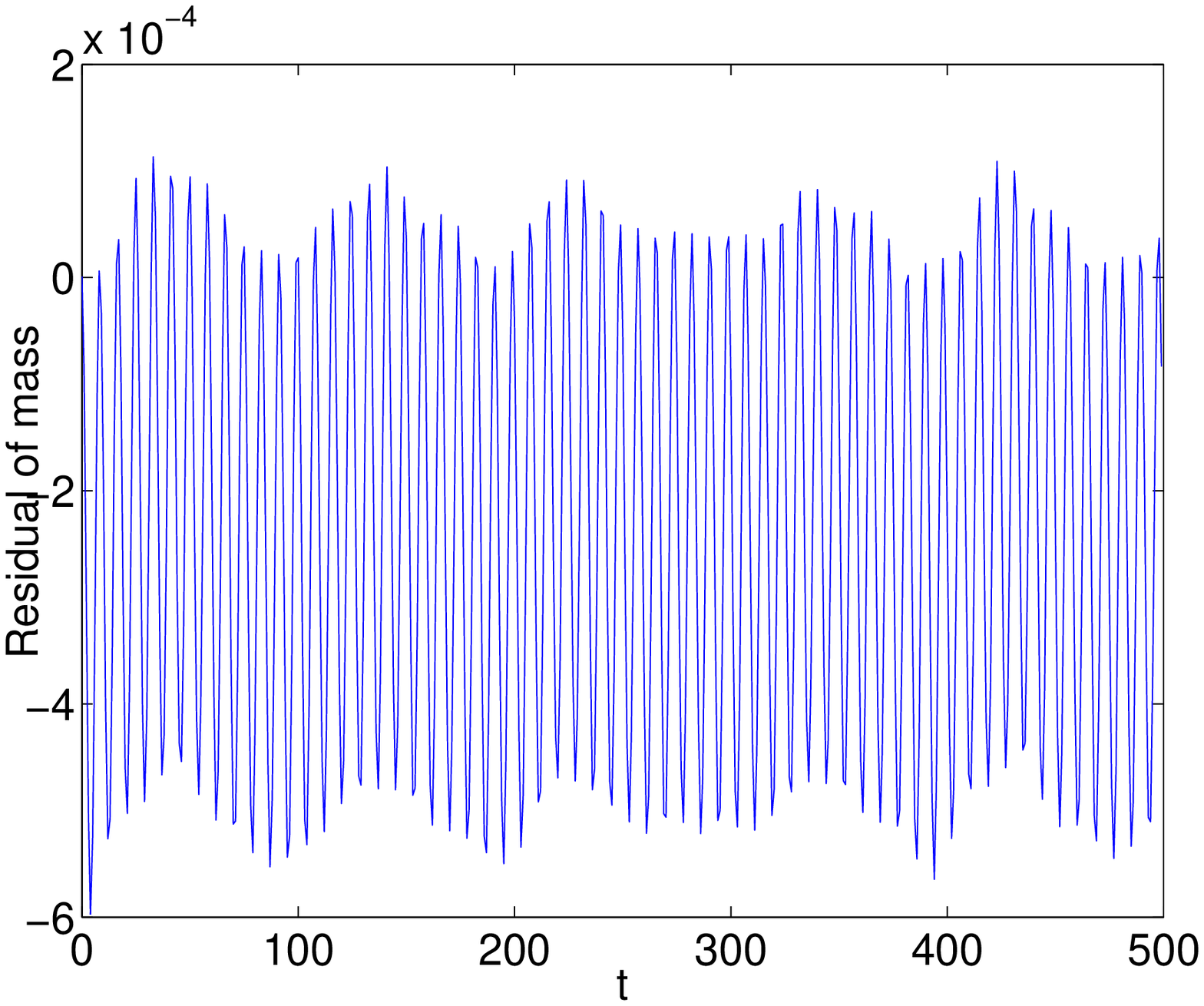,height=4.0cm,width=6.5cm},
       \epsfig{file=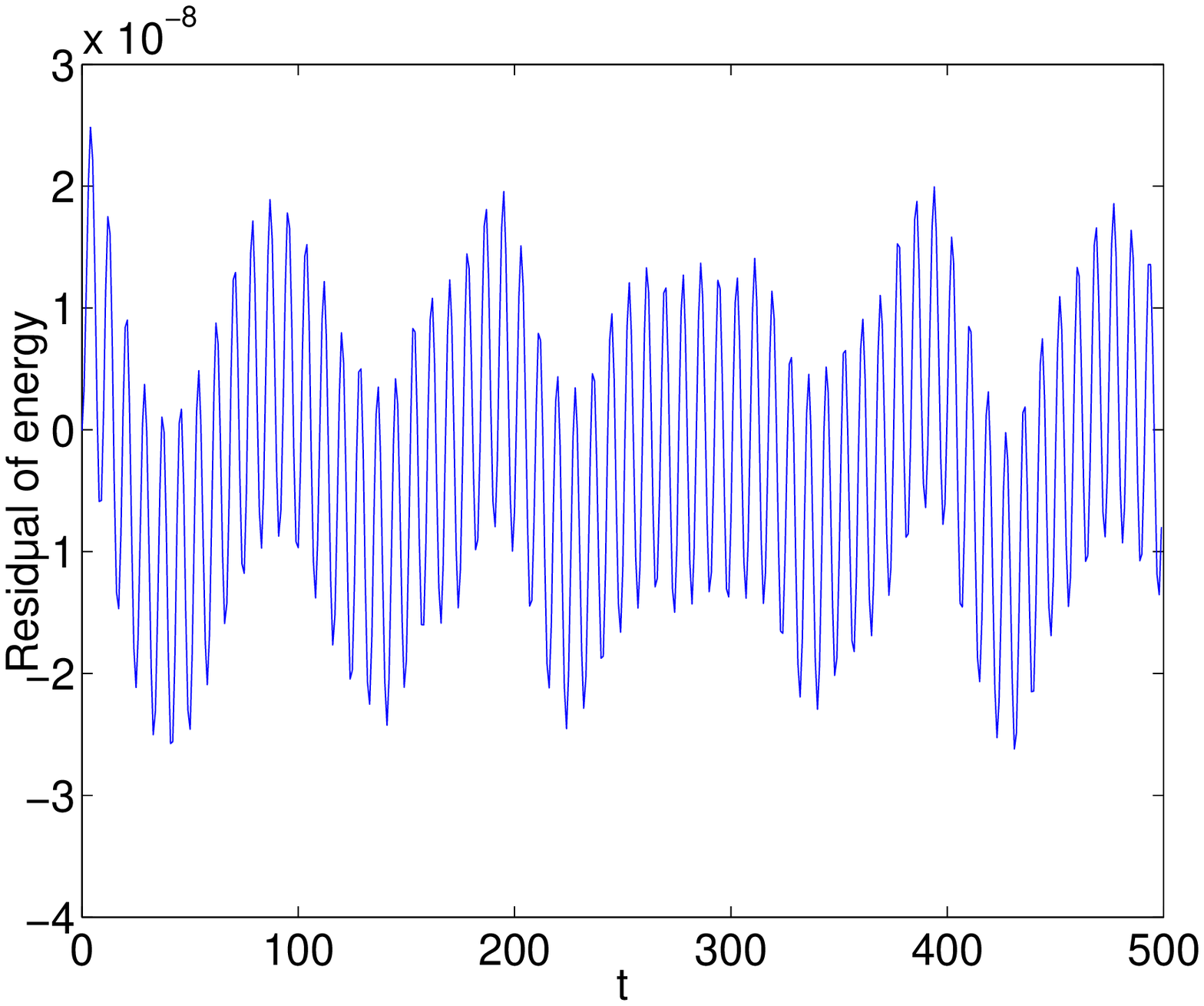,height=4.0cm,width=6.5cm}\\
       \epsfig{file=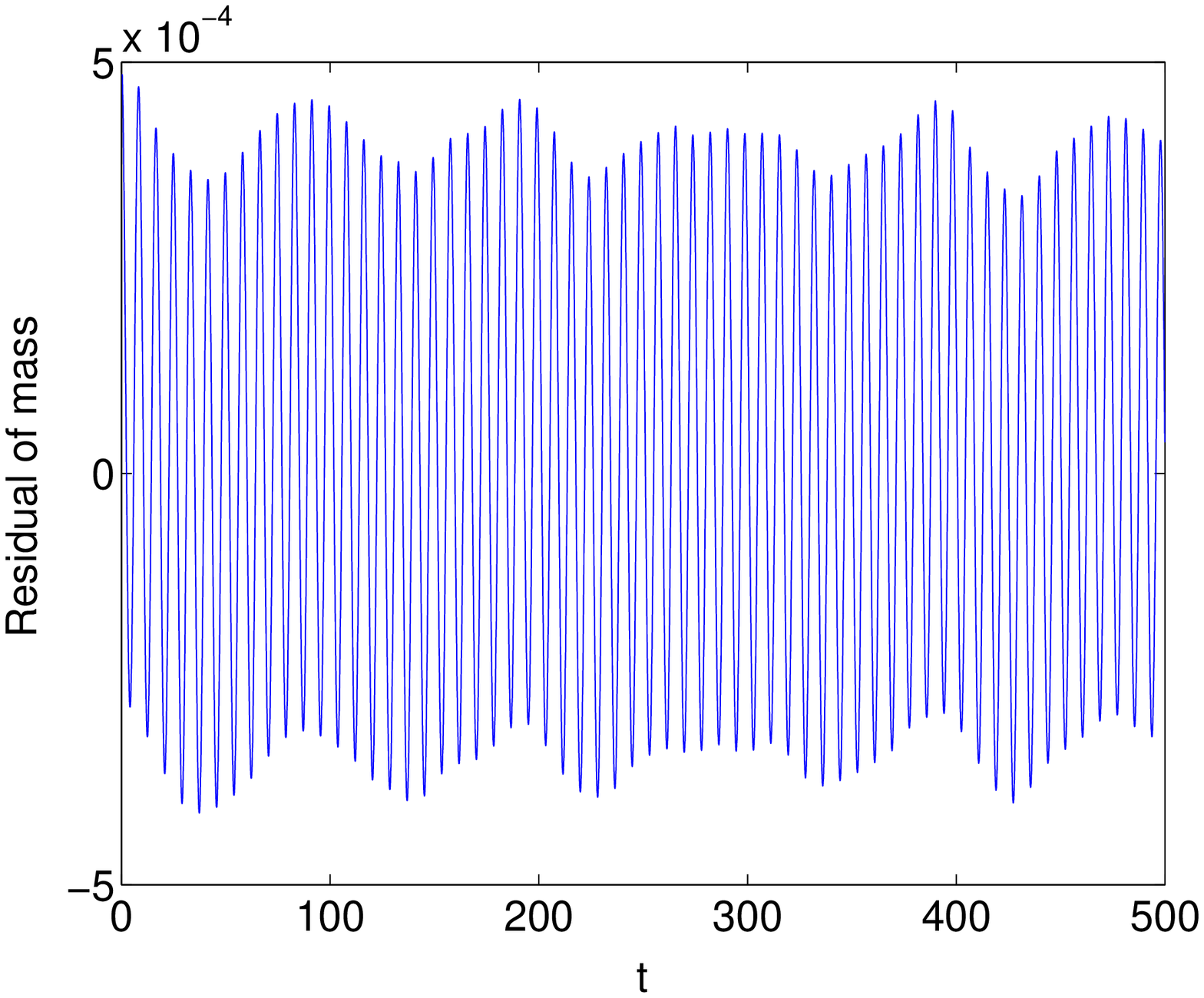,height=4.0cm,width=6.5cm},
       \epsfig{file=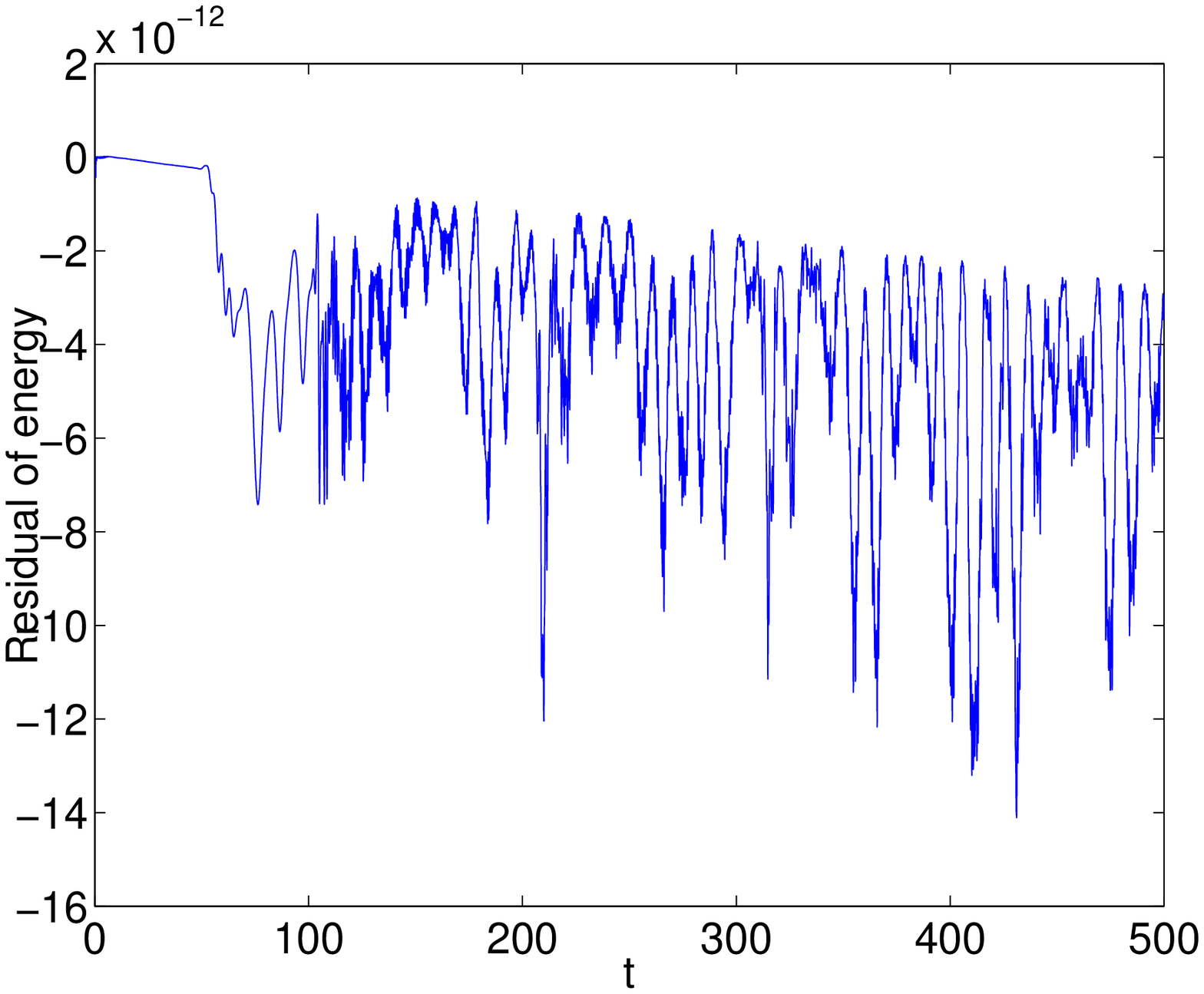,height=4.0cm,width=6.5cm}\\
       \caption{\label{compare_hom_residual}\small Conservative comparison. Upper for scheme (\ref{Preissman}), lower for scheme (\ref{scheme_ZW}). Left for mass, right for energy.}
   \end{center}
\end{figure}

The spatial-temporal domain is partitioned by $h=0.1, \tau=0.05$.
The numerical results are reported in Figs. \ref{abs_u}-\ref{compare_hom_residual}. Fig. \ref{abs_u} is the solitary wave shape of $|u|$
and its contour by MI (\ref{Preissman}). Fig. \ref{compare_hom_error} is the error of the real part of the numerical
solution by schemes (\ref{Preissman}) and (\ref{scheme_ZW}), and Fig. \ref{compare_hom_residual} is the residuals of conservative
quantities. Here we have omitted the counterpart of Fig. \ref{abs_u} by scheme (\ref{scheme_ZW}) and the error of imaginary part because they are very similar.

\begin{figure}[h]
   \begin{center}
       \epsfig{file=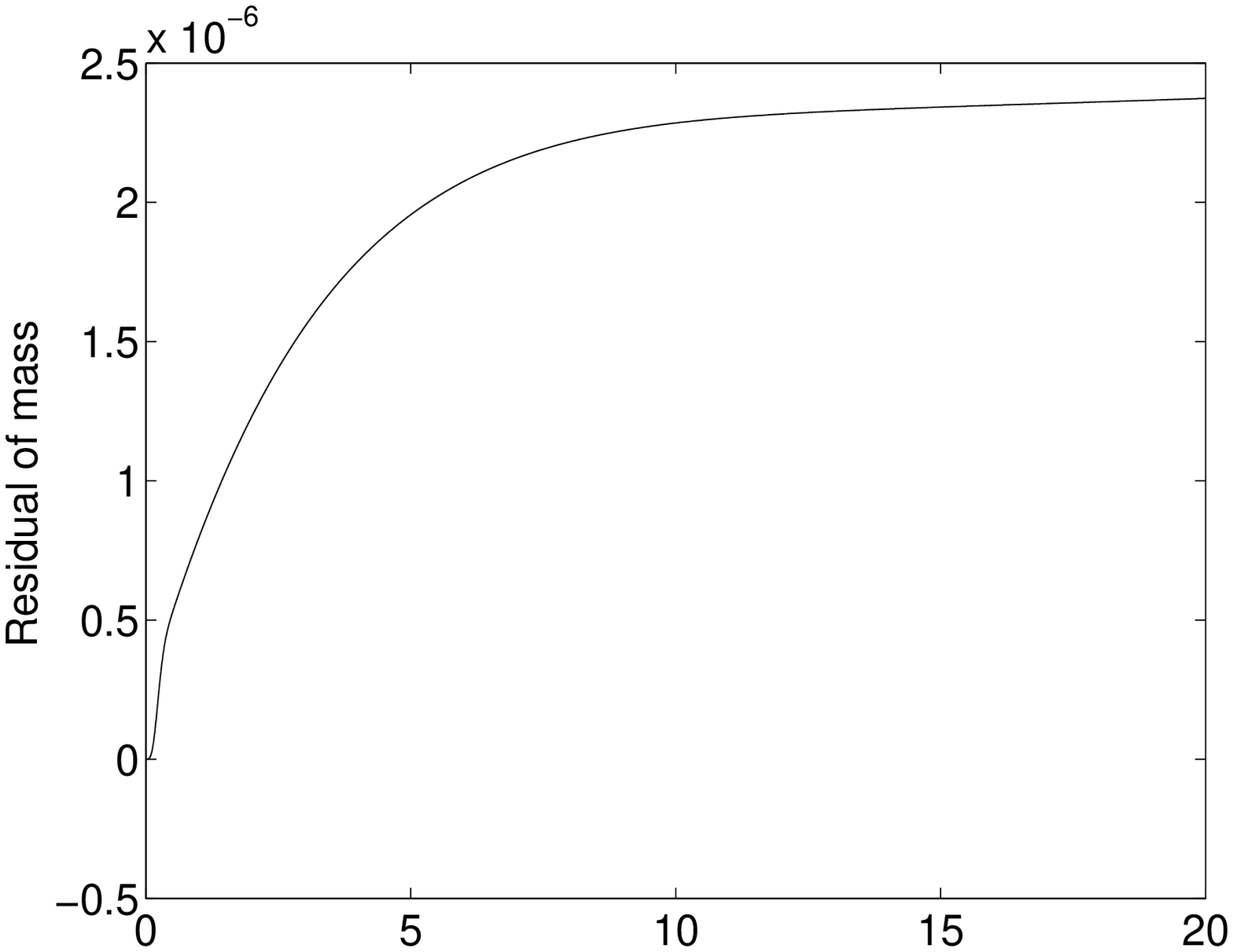,height=4.0cm,width=6.5cm},
       \epsfig{file=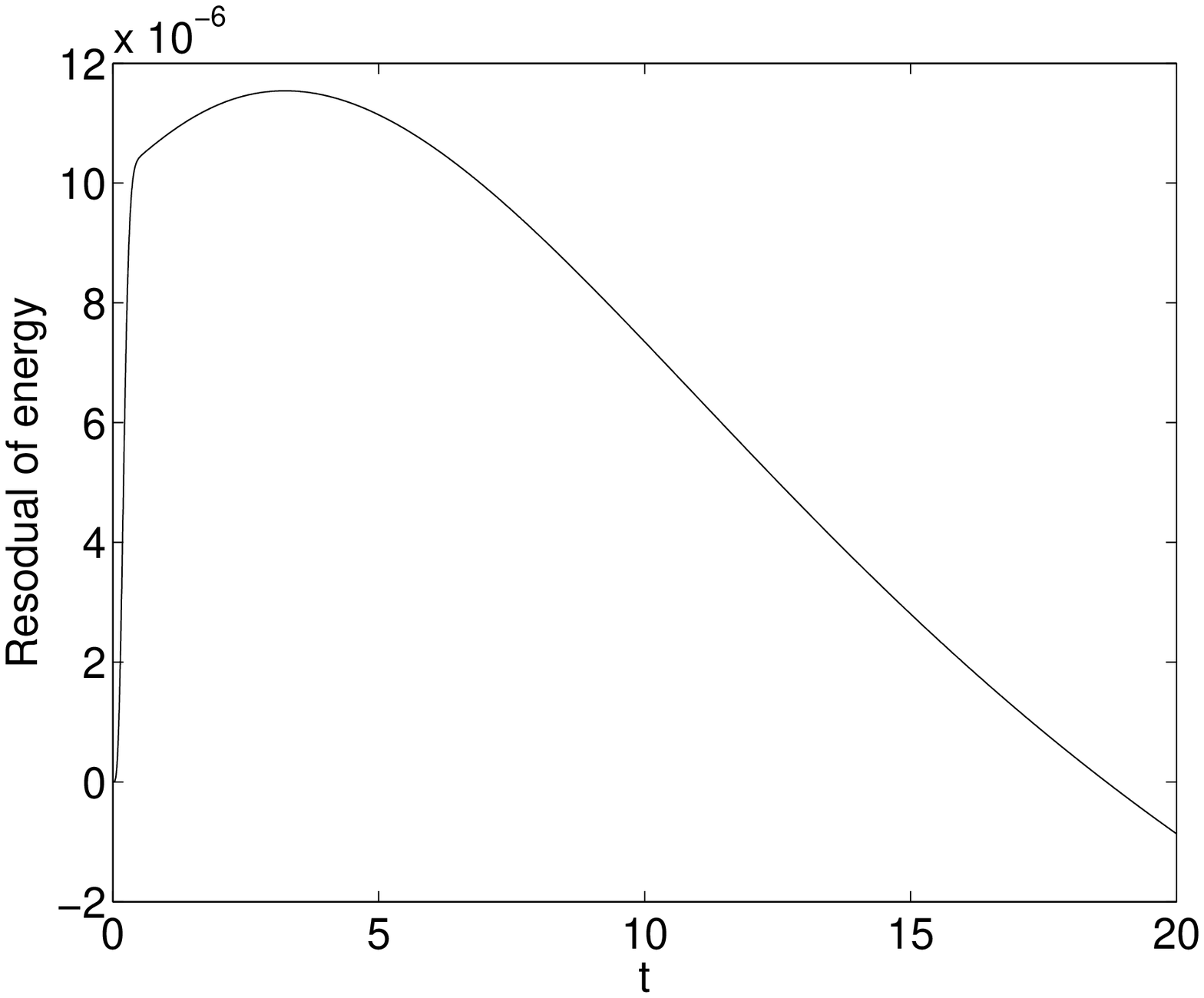,height=4.0cm,width=6.5cm}\\
       \caption{\label{Gauss_energy_mass}\small The residuals of mass (left) and energy (right).}
   \end{center}
\end{figure}

\begin{figure}[h]
   \begin{center}
       \epsfig{file=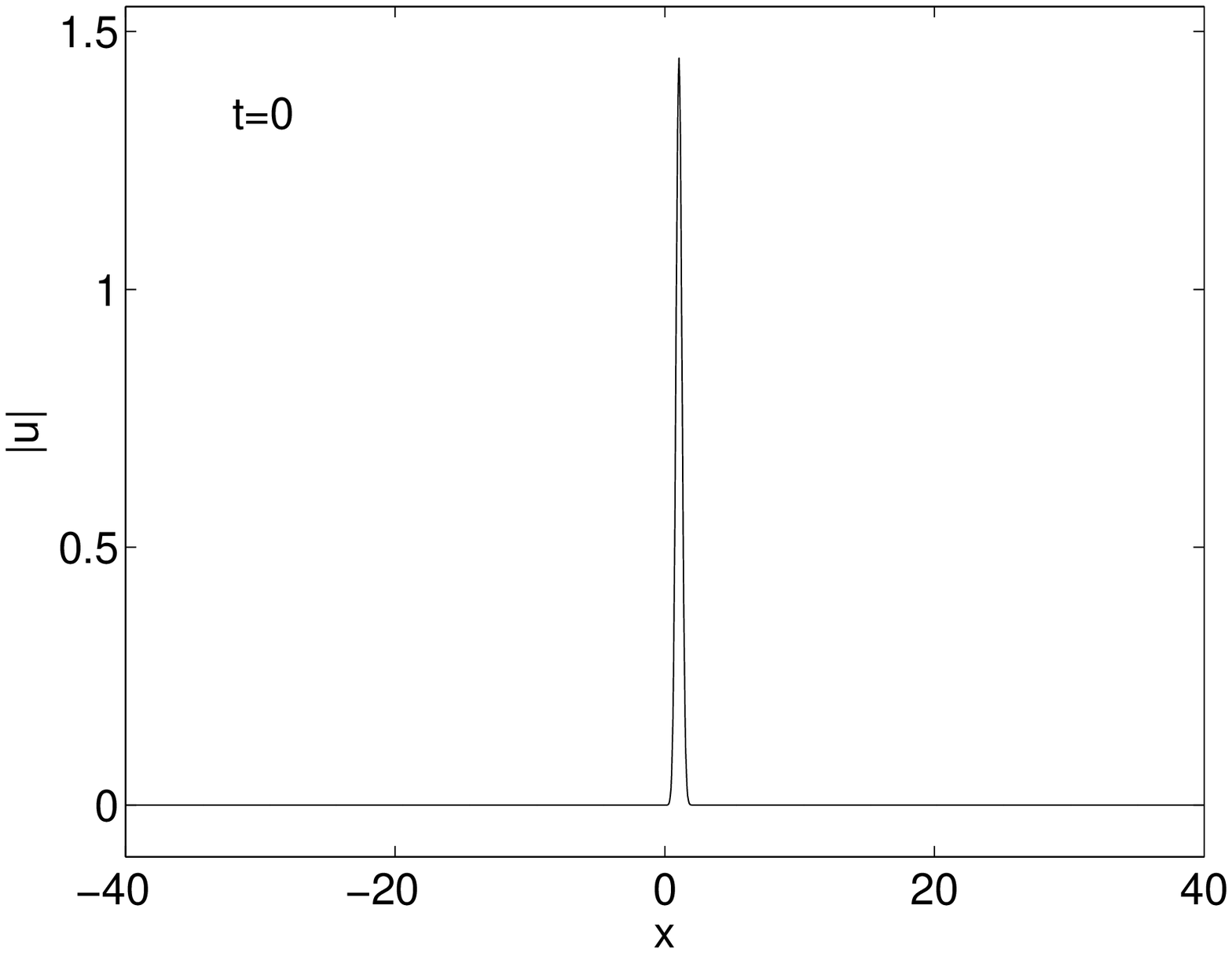,height=4.0cm,width=6.5cm},
       \epsfig{file=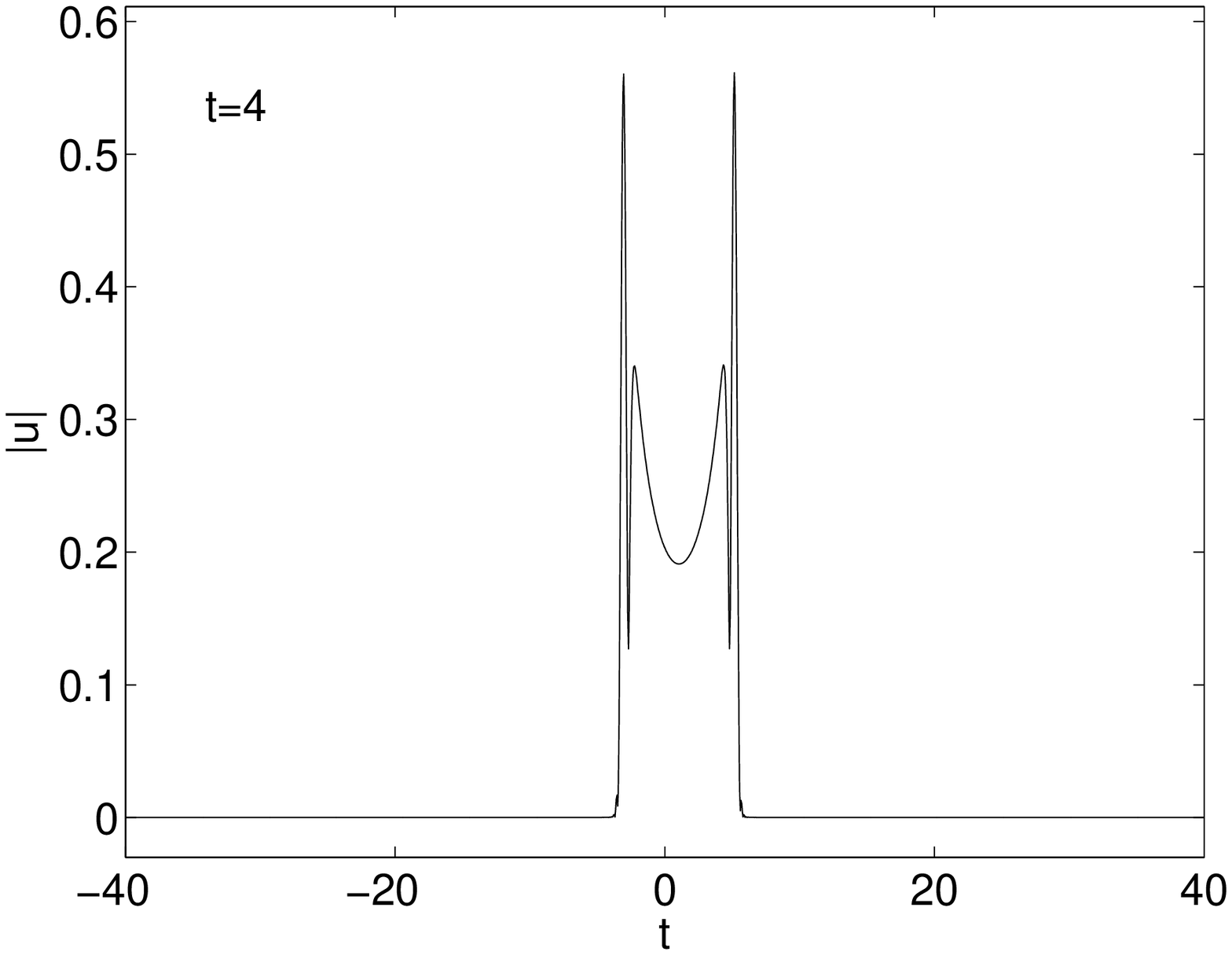,height=4.0cm,width=6.5cm}\\
       \epsfig{file=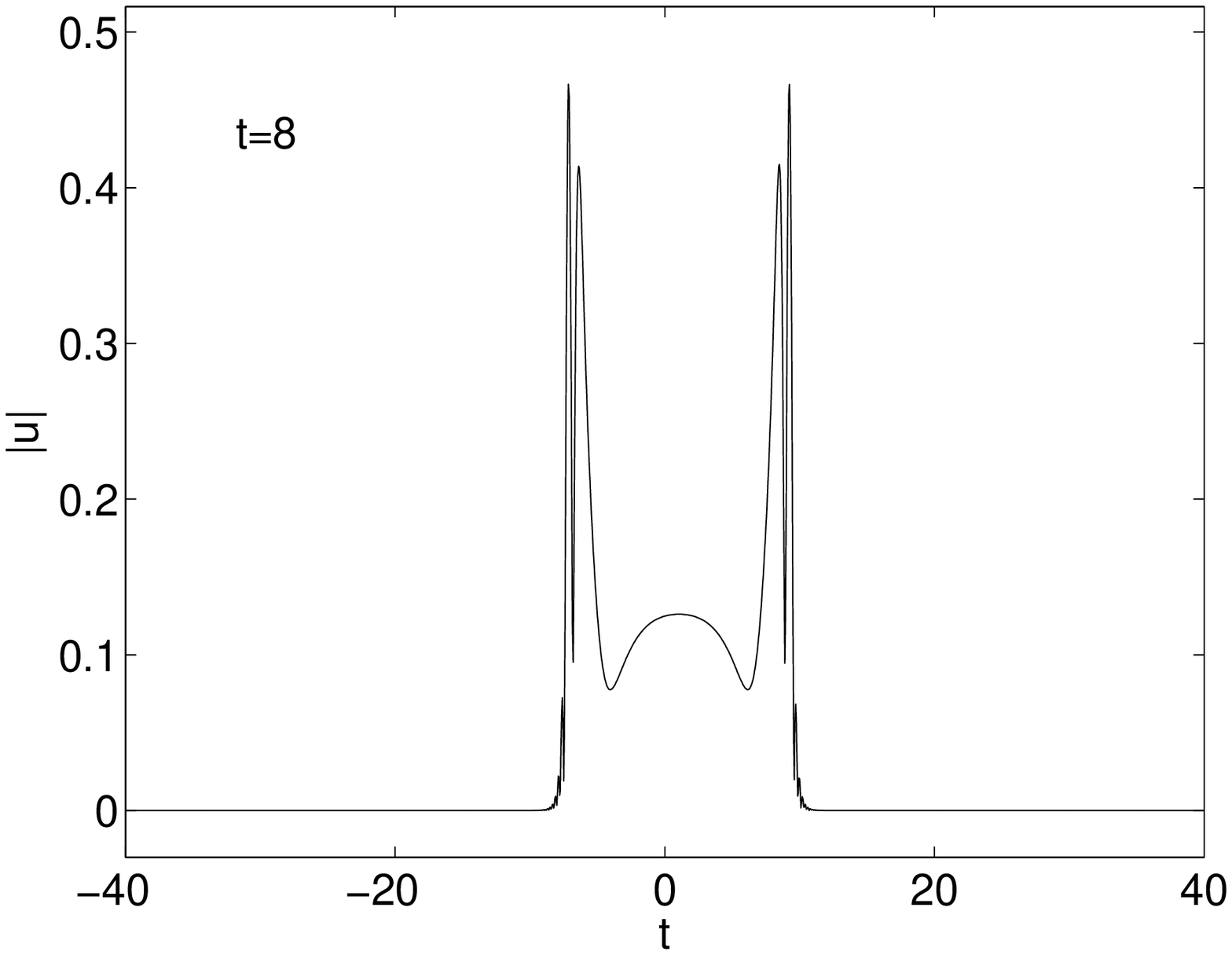,height=4.0cm,width=6.5cm},
       \epsfig{file=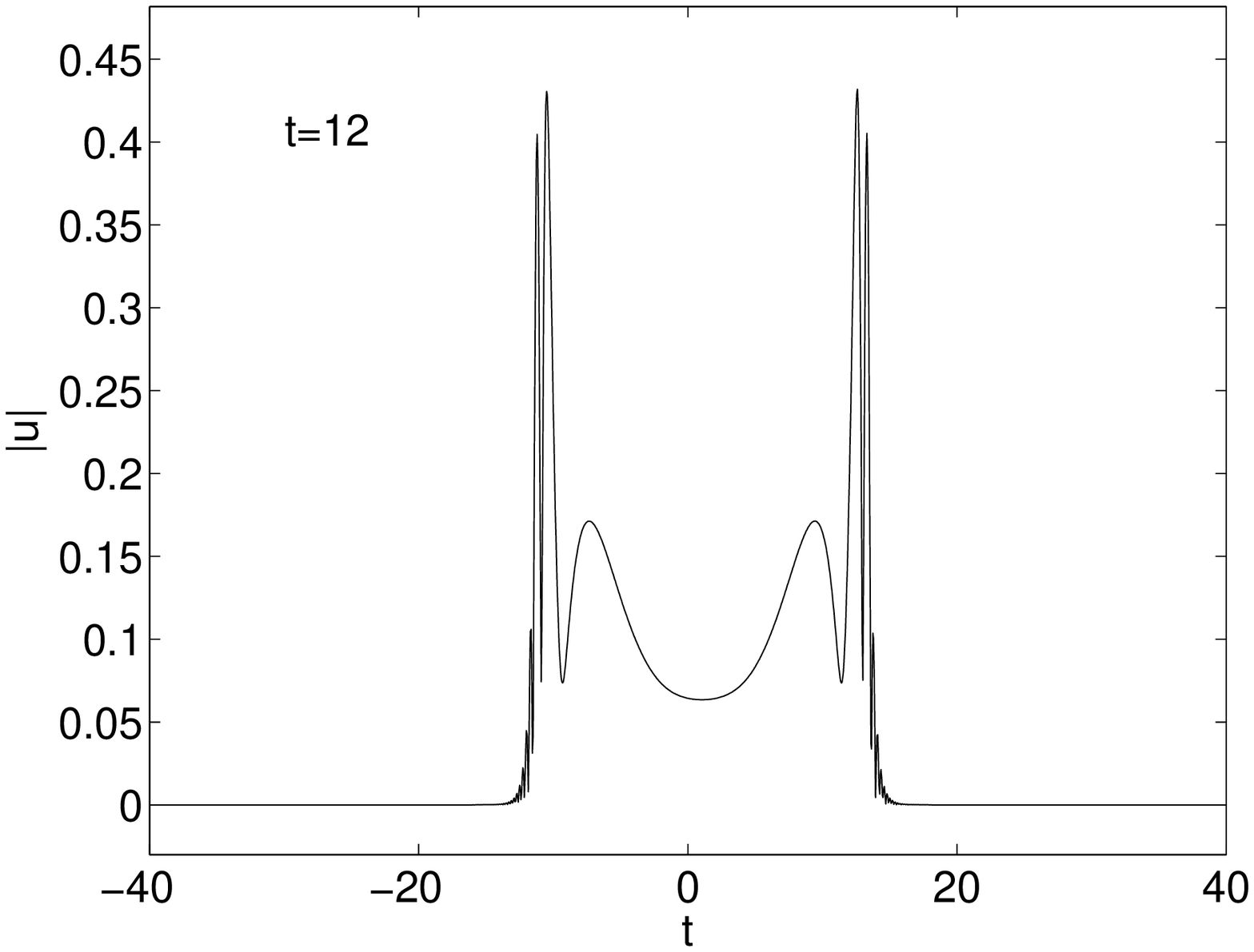,height=4.0cm,width=6.5cm}\\
       \epsfig{file=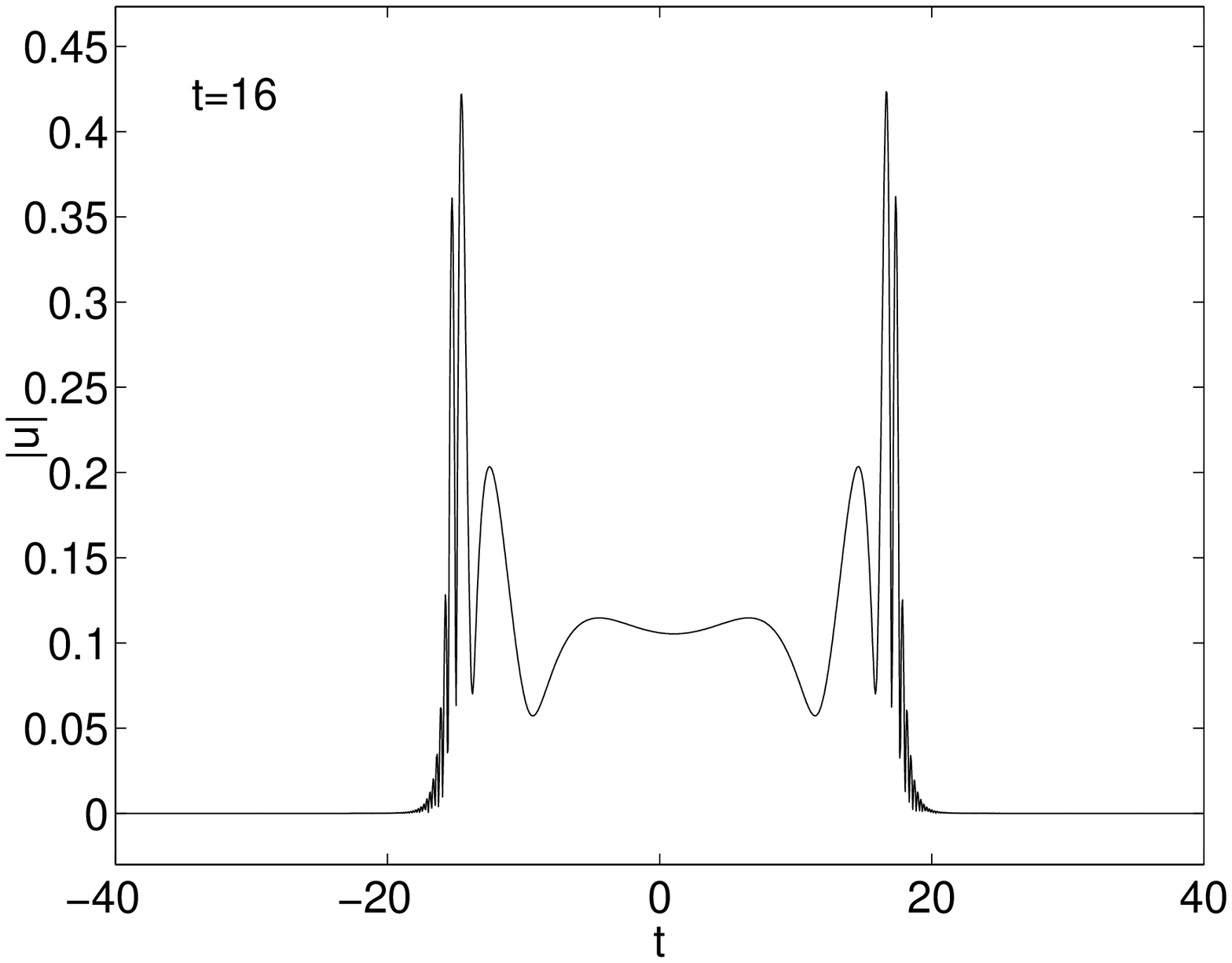,height=4.0cm,width=6.5cm},
       \epsfig{file=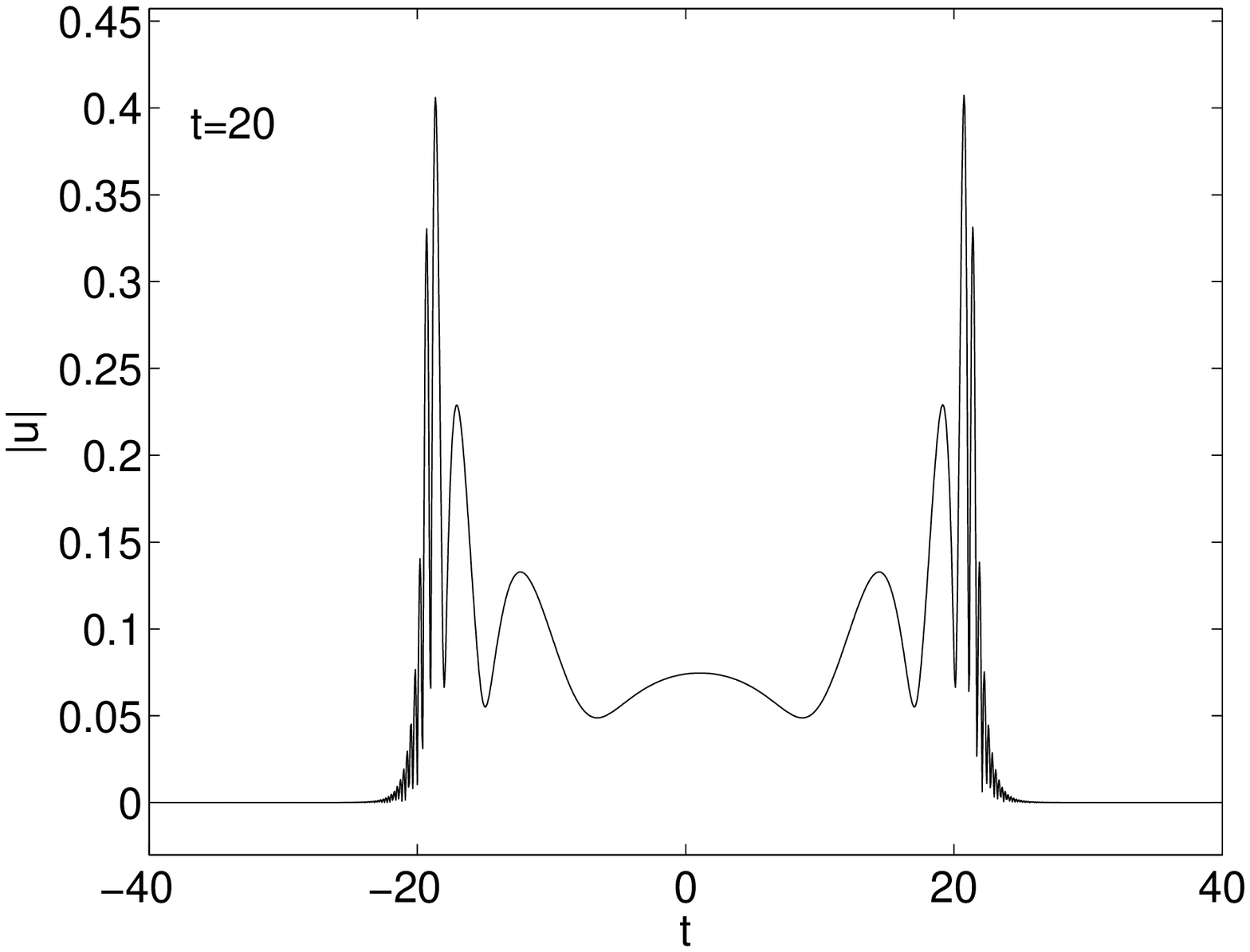,height=4.0cm,width=6.5cm}\\
       \caption{\label{Gauss_u}\small Profiles of the wave at different time stages.}
   \end{center}
\end{figure}

\begin{example}\em Finally, we consider the splitting of solitary wave\end{example}
\begin{align}\label{plane_3}
  \left\{
 \begin{array}{l}
   u_{tt}-u_{xx}+iu_t+{\left|u\right|^2}u=0,~(x,t)\in{[-40, 40]\times{(0,20]}},\\
   u(x,0)=(1+i)x\exp(-10(1-x)^2), ~~u_t(x,0)=0.
 \end{array}\right.
\end{align}
We simulate the problem by the MI (\ref{Preissman}) under the partition $h=0.05, \tau=0.02$. The residuals of mass and energy are reported in Fig. \ref{Gauss_energy_mass}, and the profiles of the numerical solution
are plotted in Fig. \ref{Gauss_u}. From the pictures, it is observed that the original wave splits into some lower waves rapidly, and more and more ripples bring up with the evolution of wave. Both
the mass and energy are nicely preserved.

\section{Conclusion}
In this study, we have used multisymplectic midpoint rule to approximate the Schr\"{o}dinger equation with
wave operator. The global conservative properties of the numerically method are investigated.
For nonlinear problem, they usually do not keep the mass and energy exactly. However, their residuals
are very small over long-term. The proposed numerical method demonstrates remarkable stable over long-term.
Through numerical illustrations, it is observed that the numerical methods is more accurate than other energy-preserving
methods.

\section*{Acknowledgement}
 This work is supported by the National Natural Science Foundation of China (Nos. 11271171, 11301234, 91130003),
 the Provincial Natural Science Foundation of Jiangxi (No. 20142BCB23009), the Foundation of Department
 of Education Jiangxi Province (No. GJJ12174), the State Key Laboratory of Scientific and Engineering
 Computing, CAS, and Jiangsu Key Lab for NSLSCS (No. 201302).

\end{document}